\documentclass[a4paper,10pt]{article}
\usepackage{geometry}
\usepackage{stmaryrd}
\usepackage{bbm}
\usepackage{amsfonts}
\usepackage{amsmath}
\usepackage{amssymb}
\usepackage{mathrsfs}
\usepackage{accents}
\usepackage{amscd}

\newtheorem{atheorem}{\bf \temp}[section]
\newtheorem{thm}[atheorem]{Theorem}
\newtheorem{cor}[atheorem]{Corollary}
\newtheorem{lem}[atheorem]{Lemma}
\newtheorem{prop}[atheorem]{Proposition}
\newtheorem{de}[atheorem]{Definition}
\newtheorem{rem}[atheorem]{Remark}

\numberwithin{equation}{section}

\title{\textbf{The perturbational stability of the Schr$\ddot{o}$dinger equation }}
\author{ Xixia Ma \footnote{Yau Mathematical Sciences Center, Tsinghua University. E-mail addresses:kfmaxixia@tsinghua.edu.cn} \ \ \ \ \ \ \
 \\}
\date{}

\begin{document}

\maketitle
\textbf{Abstract.}   By using  the Wigner transform, it is  shown that the nonlinear Schr$\ddot{\textmd{o}}$dinger equation can be described, in phase space, by a kinetic theory  similar to the Vlasov equation which is used for describing a classical collisionless plasma. In this paper we mainly show Landau damping   in the quantum sense, namely,quantum Landau damping exists for the Wigner-Poisson system. At the same time, we also prove the existence and the stability of the nonlinear  Schr$\ddot{\textmd{o}}$dinger equation under the quantum stability assumption.

\begin{center}
{\bf\large 0. \quad The origin of the problem  }
\end{center}

It is well known that classical plasma has mainly focused on regimes characterized by high temperatures and low densities, for which quantum mechanical effects have virtually on impact. however, physical systems where both plasma and quantum effects coexist do occur in nature, the most obvious example being the electron gas in an ordinary metal. At room temperature and standard metallic densities, quantum effects can no longer be ignored, so that the electron gas constitus a true quantum plasma. At the same time,  recent technological advances have it possible to envisage practical applications of plasma physics where the quantum nature of the particles plays a crucial role, such as miniaturized semiconductor devices and nanoscale objects and so on.

Now we explain  quantum effects. In fact, quantum effects can be measured by the thermal de Broglie wavelength of the particles composing the plasma $\lambda_{B}=\frac{\hbar}{mv_{T}},$ which roughly represents the spatial extension of a particle's wave function due to quantum uncertainty. For classical regimes, the de Broglie wavelength is so small that particles can be considered as pointlike, therefore there is no overlapping of the wave functions and no quantum interference. On the  other hand, it is well known from the statistical mechanics of  ordinary gas that quantum effects become important when the temperature is lower than the so-called Fermi temperature $T_{F}.$

The most fundamental model for the quantum $N-$body problem is the Schr$\ddot{o}$dinger equation for the $N-$particle wave function $\psi(x_{1},x_{2},\ldots\ldots,x_{N},t).$  We make a drastic but  useful and to some extent plausible simplification that can be achieved by neglecting  two-body (and higher order) correlations. This amounts to assume that the $N-$body wave function can be factored into the product of N one-body functions:
\begin{align}
\psi(x_{1},x_{2},\ldots,x_{N},t)=\psi_{1}(x_{1},t)\psi_{2}(x_{2},t)\ldots\psi_{N}(x_{N},t).
\end{align}
 The set of $ N $ one-body wave functions is known as a quantum mixture (or quantum mixed state) and is usually represented by a density matrix
\begin{align}
\rho(x,y,t)=\sum^{N}_{\alpha=1}p_{\alpha}\psi_{\alpha}(x,t)\psi^{\ast}_{\alpha}(y,t),
\end{align}
where, for clarity, we have assumed the same normalization $\int|\psi_{\alpha}|^{2}dx=1$ for all wave functions and then introduced the occupation probabilities $p_{\alpha}.$
And these obey $ N$ independent Schr$\ddot{o}$dinger equations, coupled through Poisson's equation
\begin{align}
i\hbar\frac{\partial\psi_{\alpha}}{\partial t}=-\frac{\hbar^{2}}{2m}\frac{\partial^{2}\psi_{\alpha}}{\partial x^{2}}-e\phi\psi_{\alpha},\quad\quad \alpha=1,\ldots,N,
\end{align}
\begin{align}
\frac{\partial^{2}\phi}{\partial x^{2}}=\frac{e}{\varepsilon_{0}}\bigg(\sum^{N}_{\alpha=1}p_{\alpha}|\psi_{\alpha}|^{2}-n_{0}\bigg),
\end{align}
where $\hbar$ is the Plank constant, m is the mass of the particles and $e$ is the charge of the particles.


 Under the above background, we try to describe the quantum Landau damping, using analytical theory and simulations, and  assess the impact of quantum effects on the nonlinear wave-particle interactions in quantum electrodynamic field. The start-point of this problem comes from the basic physics fact: wave-particle phenomenon in nature. In classical plasma physics, the notion of wave-particle interactions, the couplings between collective and individual particles' behaviors, is fundamental to our comprehension of plasma phenomenon. In particularly, Landau damping factor is obtained from actual  physical process based on wave-particle interaction.
In this paper we try to study  the behavior of quantum plasma under the quantum electric field,  and if  the asymptotic behavior of the corresponding Wigner function is exponential decay,  we call that  as  quantum Landau damping. In fact, quantum Landau damping is the phenomenon that describes the probability of the distribution of  the single particle's motion.  And we will prove it   from the rigorous  mathematical sense.

From now on, we restrict the motion of the charged particle  in a periodic box of length 1, that is, a torus $\mathbb{T}^{3}$ of length 1, containing a uniform motionless background of irons.


We give a detailed kinetic  description of the motion model (0.3)-(0.4) of quantum plasma. And note that $\sum^{N}_{\alpha=1}\int_{\mathbb{T}^{3}_{x}}|\psi_{\alpha}|^{2}dx$ is conserved by the motion.  Indeed, multiplying (0.3) by the complex conjugate $\bar{\psi}_{\alpha}$ of the wave function $\psi_{\alpha},$ integrating over $\mathbb{T}^{3}_{x}$ and taking imaginary parts gives:
$\frac{d}{dt}\int_{\mathbb{T}^{3}_{x}}|\psi_{\alpha}|^{2}dx=0.$ Then the normalization of $\psi$ is taken as $\sum^{N}_{\alpha=1}\int_{\mathbb{T}^{3}_{x}}|\psi_{\alpha}|^{2}dx=\sum^{N}_{\alpha=1}\int_{\mathbb{T}^{3}_{x}}\bar{\psi}_{\alpha}\psi_{\alpha }dx=1.$




Denote $$H_{x}\doteq \frac{\hbar^{2}}{2m}\Delta_{x}+e\phi,$$ which stands for the Hamiltonian acting on the $x.$

Differentiating (0.2) with respect to $t$, and then  using  (0.3) gives the evolution for the density matrix $\rho(t,s,r):$
\begin{align}
i\hbar\partial_{t}\rho(t,s,r)=(H_{s}-H_{r})\rho(t,s,r).
\end{align}
And the equation (0.7) reads explicitly
\begin{align}
i\hbar\partial_{t}\rho=-\frac{\hbar^{2}}{2m}(\Delta_{s}-\Delta_{r})\rho-e(\phi(s,t)-\phi(r,t))\rho.
\end{align}

Then we reformulate the quantum equations of the  motion in kinetic form. First, we introduce the change of coordinates
\begin{align}
r=x+\frac{\hbar}{2m}\eta,\quad s=x-\frac{\hbar}{2m}\eta
\end{align}
and set
\begin{align}
u_{\hbar}(t,x,\eta)=\rho(t,x+\frac{\hbar}{2m}\eta,x-\frac{\hbar}{2m}\eta).
\end{align}

 Second, we  consider Fourier transforms of functions which depend on $ \eta.$ Since $(\hbar/2m)\eta$ has the dimension of length, we conclude that $\eta$ has the dimension of inverse velocity, and thus, the dual variable of $\eta$ has the dimension of velocity.
Therefore, the  function $w_{\hbar},$ which corresponds to the wave function $\psi$ (or, equivalently, to the density matrix $\rho$ given by (0.2)) is defined as the inverse Fourier transform of $ u_{\hbar}$ with respect to $\eta:w_{\hbar}:=\mathcal{F}^{-1}u_{\hbar}$ or, explicitly:
$$w_{\hbar}(t,x,v)=\frac{1}{(2\pi)^{3}}\int_{\mathbb{R}^{3}}\rho(t,x+\frac{\hbar}{2m}\eta,x-\frac{\hbar}{2m}\eta)e^{i\eta\cdot v}d\eta,$$
that is,
\begin{align}
w_{\hbar}(t,x,v)=\frac{1}{(2\pi)^{3}}\sum^{N}_{\alpha=1}p_{\alpha}\int_{\mathbb{R}^{3}}\overline{\psi_{\alpha}(t,x+\frac{\hbar}{2m}\eta)}
\psi_{\alpha}(t,x-\frac{\hbar}{2m}\eta)e^{iv\cdot\eta}d\eta.
\end{align}
It was introduced by E.Wigner  in 1932 and,called the Wigner transform, which is a useful tool to express quantum mechanics in a phase space formalism. As above, a generic quantum mixed state can be described by $ N $ single particle wave function $\psi_{\alpha}(x,t)$ each characterized by a probability $p_{\alpha}$ satisfying $\sum^{N}_{\alpha=1}p_{\alpha}=1.$   We shall see below, its construction constitutes a major break-through in the quest for a kinetic formulation of quantum transport. The Wigner function is a function of the phase space variables $(x,v)$ and time, which, in terms of the single-particle wave functions.

 We note that the Wigner function behaves as the classical particle distribution. Using (0.5), we immediately derive that the mean value of the Wigner function $w_{\hbar}$ with respect to the velocity $v$ is the quantum electron ensemble position density
\begin{align}
n_{\hbar}(t,x)=\int_{\mathbb{R}^{3}}w_{\hbar}(t,x,v)dv,
\end{align}
since $u_{\hbar}(t,x,\eta=0)=\mathcal{F}w_{\hbar}(t,x,\eta=0)=\rho(t,x,x)$ holds. However, it must be stressed  that  the Wigner function does not necessarily stay nonnegative in its evolution process. Unlike in the classical case, it can therefore not be interpreted as a probability density. In the literature it is often referred to as `quasi-distribution' of particles.

The evolution equation for the Wigner functions is obtained by transforming  the Heisenberg equation (0.7) for the density matrix $\rho$ to the $(x,\eta)-$coordinates given by (0.8):
\begin{align}
\partial_{t}u_{\hbar}+\textmd{idiv}_{\eta}(\nabla_{x}u_{\hbar})+ie\frac{\phi(t,x+\frac{\hbar}{2m}\eta)-\phi(t,x+\frac{\hbar}{2m}\eta)}{\hbar}u_{\hbar}=0,
\end{align}
and by the  Fourier transformation
\begin{align}
\partial_{t}w_{\hbar}+v\cdot\nabla_{x}w_{\hbar}+\frac{e}{m}\Theta_{\hbar}[\phi]w_{\hbar}=0,
\end{align}
where the operator $\Theta_{\hbar}[\phi]$ is defined by
$$(\Theta_{\hbar}[\phi]w_{\hbar})(t,x,v)$$
\begin{align}
=\frac{im}{(2\pi)^{3}}\int_{\mathbb{R}^{3}}\int_{\mathbb{R}^{3}}\frac{\phi(t,x+\frac{\hbar}{2m}\eta)-\phi(t,x-\frac{\hbar}{2m}\eta)}{\hbar}w_{\hbar}(t,x,v')e^{i(v-v')\cdot\eta}
dv'd\eta,
\end{align}
where $\phi(x,t)$ is the self-consistent electrostatic potential,
and the corresponding initial data $w_{\hbar I}$ is \begin{align}
w_{\hbar I}(x,v)=\frac{1}{(2\pi)^{3}}\int_{\mathbb{R}^{3}}\overline{\psi_{I}(x+\frac{\hbar}{2m}\eta)}\psi_{I}(x-\frac{\hbar}{2m}\eta)e^{iv\cdot\eta}d\eta
\end{align}

Developing the integral term in (0.12) up to order $O(\hbar^{2}),$ we obtain
\begin{align}
\frac{\partial w_{\hbar}}{\partial t}+v\cdot\nabla_{x}w_{\hbar}-\frac{e}{m}\nabla_{x}\phi\cdot\nabla_{v}w_{\hbar}=\frac{e\hbar^{2}}{24m^{3}}\frac{\partial^{3}\phi}{\partial x^{3}}\frac{\partial^{3}w_{\hbar}}{\partial v^{3}}+O(\hbar^{4}).
\end{align}
The Vlasov equation is thus recovered in the formal semiclassical limit $\hbar\rightarrow0.$ We stress that rigorous asymptotic results are much harder to obtain and generally involve weak convergence. At present, there are some works on it, such as [1,24,].

The Wigner equation must be coupled to the Poisson's equation for the electric potential
\begin{align}
\frac{\partial^{2}\phi}{\partial x^{2}}=\frac{e}{\varepsilon_{0}}\bigg(\int w_{\hbar}dv-n_{0}\bigg)
\end{align}
where we have assumed that the ions form a motionless neutralizing background with uniform density $n_{0}.$

 Comparing with the Vlasov equation for a collisionless particle plasma, this scheme's advantage is basically connected to the well-known property of quantum diffraction to prevent the too steep gradients in phase space appearing in usual codes, by simply reorganizing the finer and finer structure at the level of the elementary cell $\hbar.$ The above Wigner-Poisson(WP) system has been extensively used in the study of quantum transport. Exact analytical results can be obtained by linearizing (0.14) and (0.18) around a spatially homogeneous equilibrium given by $w^{0}_{\hbar}(x,v),$ where $w^{0}_{\hbar }$  is defined as follows: \begin{align}
w^{0}_{\hbar}(x,v)=\frac{1}{(2\pi)^{3}}\int_{\mathbb{R}^{3}}\overline{\Psi(x+\frac{\hbar}{2m}\eta)}\Psi(x-\frac{\hbar}{2m}\eta)e^{iv\cdot\eta}d\eta.
\end{align}

First, We recall the development of the classical Landau damping.  It is well known that the existence of a damping mechanism by which plasma particles
absorb wave energy  was found by L.D.Landau at the linear level, under the condition that the plasma is not cold and the velocity distribution is of finite extent. and later we call it Landau damping. A ground-breaking work for Landau damping was made by Mouhot and Villani in the  nonlinear case.
 They gave the first and rigours proof of nonlinear Landau damping under some assumption of the electric field. Later Masimudi and so on gave a simpler proof.
 Recently, based on the physical viewpoint, we remove the stability condition introduced by  Mouhot and Villani, and
 characterized the corresponding  stability condition from the energy perspective.
  $$ \textit{Based } \textit{on} \textit{ the }\textsl{above } \textit{facts } \textit{and } \textit{observation } \textit{in  }\textit{classical } \textit{plasma, } \textit{we}\textit{ claim }\textit{that, }\textit{in }  \textit{nonlinear } \textit{quantum }$$
  $$\textit{electrodynamical }\textit{cases, } \textit{ quantum }\textit{Landau }\textit{damping } \textit{similar }\textit{to }\textit{classical } \textit{Landau }\textit{damping }\textit{still }\textit{occurs, }$$
  $$\textit{although }\textit{of }\textit{course, }\textit{the }\textit{numerical }\textit{value }\textit{of }\textit{the }\textit{damping }\textit{rate }\textit{will }\textit{depend }\textit{on }\textit{$\hbar.$}
  \quad\quad\quad\quad\quad\quad\quad\quad\quad\quad$$

  As a valid evidence of the above claim, there have been lots of experiment's results which can be  found in [7,11,12,27,28,30].

  However, there are still many differences between the classical case and the quantum case from the conditions of occurrence or technical methods. The main difference is that  the density function $n=|\psi|^{2}$ obtained from the wave function in the SP equations does not recover the classical density field $\rho=\int f dv$ from the classical VP equations.  Indeed, we consider a distribution which is the superposition of two Gaussian waves traveling with opposite velocities in the quantum case; however, in the classical case, the superposition is just a Gaussian. Then the density in the quantum case always exhibits order unity differences from the classical solution, and the potential(obtained from the density via the Poisson equation) and the force field (gradient of the potential)will also oscillate whenever in linear case or in nonlinear case.




 In this paper in analogy with the classical case,  we can regard the Wigner function as a `distribution' function of a single wave. Here $v$ is  the particle's velocity and  $v_{T}$ is the thermal  velocity.

 Now  we recall Landau damping through  gathering lots of physical literature and results of mathematical articles, mainly in classical plasma, since the literature on the stability  is  very scarce from the mathematical viewpoint as far as the case in quantum plasma.
   Except Landau's work in linear case, many works from mathematical aspects found in [6,16,23,32]
 gave rigorous proofs under different assumptions. Later,
 In recent years, there are also lots of results on the stability in other setting such as those in [2,5,10,13,14,17,18,33,34].

In this paper  we will prove the main results in the following.
\begin{thm} Let $w_{\hbar}^{0}:\mathbb{R}^{3}\rightarrow\mathbb{R}_{+}$ be an analytic velocity profile,
 and let $\phi(t,x)=W(x)\ast n(t,x)$ and $W(x):\mathbb{T}^{3}\rightarrow\mathbb{R}$ satisfy $$ \widehat{W}(0)=0,\quad
  |\widehat{W}(k)|\leq\frac{1}{1+|k|^{\gamma}},\gamma\geq1.$$ Further assume that, for some constants $\bar{\lambda}_{0}>\lambda_{0}>0,\bar{\mu}_{0}>\mu_{0}>0,$ the functions $w_{\hbar I},w^{0}_{\hbar}$ defined in (0.16)-(0.17) satisfy the following conditions:
 \begin{align}
 \sup_{\eta\in\mathbb{R}^{3}}e^{2\pi\lambda_{0}|\eta|}|\tilde{w}_{\hbar}^{0}(\eta)|\leq C_{0},\quad \sum_{n\in\mathbb{N}_{0}^{3}}\frac{\lambda_{0}^{n}}{n!}\|\nabla^{n}_{v}w_{\hbar}^{0}\|_{L_{dv}^{1}}\leq C_{0}<\infty.
 \end{align}
  and there is $\varepsilon_{0}=\varepsilon_{0}(\lambda_{0},\mu_{0},\beta,\gamma,\lambda'_{0},\mu'_{0})$ such that
 \begin{align}
\sup_{k\in\mathbb{Z}^{3},\eta\in\mathbb{R}^{3}}e^{2\pi\lambda_{0}|\eta|}e^{2\pi\mu_{0}|k|}|w_{\hbar}^{0}-w_{\hbar I}|+\int_{\mathbb{T}^{3}}\int_{\mathbb{R}^{3}}
|w_{\hbar}^{0}-w_{\hbar I}|e^{\beta|v|}dvdx\leq\varepsilon
\end{align}
 holds for any $\beta>0.$

At the same time, we also assume that for any $l\in\mathbb{Z}^{3},$  the following stability condition holds:
 set $$\delta w^{0}_{\hbar}(t,x,v)=\int_{\mathbb{T}^{3}} e^{-ix\cdot k} \frac{w^{0}_{\hbar}(k,v-\frac{\hbar}{2m}l)-w^{0}_{\hbar}(k,v+\frac{\hbar}{2m}l)}{\frac{\hbar}{m}l}dk.$$
\begin{itemize}


             \item[] $\mathbf{Stability}$ $ \mathbf{condition}:$  there is some  constant $\delta_{0}>0,$  \begin{align}
\bigg\|\int_{\mathbb{T}^{3}}\delta w^{0}_{\hbar}(x,\cdot)dx\bigg\|_{\mathcal{C}^{\bar{\lambda}(1+b);1}}\leq \delta_{0},\quad\quad
\sup_{0<\tau\leq t}\|\delta w^{0}_{\hbar}\|_{\mathcal{Z}^{\bar{\lambda}(1+b),\bar{\mu};1}_{\tau-\frac{bt}{1+b}}}\leq\delta_{0},
\end{align}
         \end{itemize}
         Then 
for any fixed $\eta,k,\forall r\in\mathbb{N},$ as $|t|\rightarrow\infty,$ we have
         \begin{align}
         &|\hat{w}_{\hbar}(t,k,\eta)-\hat{w}_{\hbar I}(k,\eta)|\leq  e^{-\lambda'_{0}|\eta+kt|},\quad\|n_{\hbar}(t,\cdot)-n_{\hbar I}\|_{C^{r}(\mathbb{T}^{3})}
         =O(e^{-2\pi\lambda'_{0}|t|}),\notag\\
          &\|\phi(t,\cdot)\|_{C^{r}(\mathbb{T}^{3})}=O(e^{-2\pi\lambda'_{0}|t|}),\notag\\
         \end{align}
         where $n_{\hbar I}=\int_{\mathbb{T}^{3}}\int_{\mathbb{R}^{3}}w_{\hbar I}(x,v)dvdx,$ for any $ 0<\lambda'_{0}<\lambda_{0} .$
         \end{thm}

          \begin{rem}
          Now we simply analyze the differences both  the above stability and the classical stability. First,  recall the stability condition in  the classical plasma: $$ \textmd{for } \textmd{any }\textmd{velocity }  v\in\mathbb{R}^{3},
            \textmd{ there } \textmd{exists } \textmd{some }\textmd{ positive }\textmd{ constant } v_{T}\in\mathbb{R}\textmd{ such } \textmd{that }\textmd{if } v=\frac{\omega}{k},\textmd{ then }|v|\gg v_{T}. $$

             Comparing to the classical case, on the one hand, the quantum stability condition implies that the electron wave behaves as a quasi-particle; the increase or decrease in the wave energy by one quantum is accompanied by the absorption or emission of electron energy $\hbar\omega$ and momentum $\hbar k.$  Now we give a simple calculation: for the momentum p,
$$p=p\cdot\frac{k}{k}=_{q}= m \frac{\omega}{k}\mp\frac{\hbar k}{2}=_{c,sc}m \frac{\omega}{k};$$
for the energy changes e,
$$e=\pm\hbar\omega=_{q}\pm\hbar k\cdot\frac{p}{m}+\frac{\hbar^{2}k}{2m}=_{c,sc}\pm\hbar k\cdot\frac{p}{m}.$$
From  the first equation, it reveals an important effects of the wave-like,diffractive nature of electrons. That is, in the quantum case, the resonant momentum condition is shifted by $\mp\frac{\hbar k}{2}$ with respect to the classical and semi-classical conditions, this in turn affects the relation of the Landau damping rate. The second equality implies that the energy transfer differs by the reoil energy $E_{rec}(k)=\frac{\hbar^{2}k^{2}}{2m}$ associated with the momentum transfer $\pm\hbar k.$

On the other hand,
 in the classical case, the superposition is just a Gaussian. However,  we  have to consider a distribution which is the superposition of two Gaussian waves traveling with opposite velocities in the quantum case, which lead the resonance occurring even in the linear case.

         \end{rem}

           \begin{center}
\item\section{Physical deduction of quantum Landau damping}
\end{center}

In this section our goal is to remind the reader that the classical damping originates from the singularity appearing in the following dispersion relation (1.6) at the point $v=\omega/k$ in velocity space. This corresponds to particles whose velocity is equal to the phase velocity of the wave $\omega/k$ (resonant particles). In the following we show that the argument which is originally developed for the Vlasov-Poisson system, still holds for the quantum Wigner-Poisson case, although the damping rate will depend on $\hbar.$

Accordingly the  field to be weak, we put $\rho=\rho_{0}(r_{1}-r_{2})+\delta\rho(t,r_{1},r_{2}),$ where $\rho_{0}$ is the density matrix of the unperturbed stationary and homogeneous state of the electron; the homogeneity implies that $\rho_{0}$ depends only on the difference $R=r_{1}-r_{2}.$ The density matrix $\rho_{0}(R)$ is related to the (unperturbed)electron momentum distribution function $n_{0}(p)=\mathcal{N}_{e}\int\rho_{0}(R)e^{-ip\cdot R/\hbar}d^{3}x,$ where
$\mathcal{N}_{e}$ is the total number of electrons. The function $n_{0}(p)$ is defined as the occupation numbers of quantum states of electrons with definite values of the momentum and the spin component. The number of states in an element $d^{3}p$ of momentum space is $2d^{3}p/(2\pi\hbar)^{3}.$
If omitting terms of the second order of smallness, we obtain a linear equation for the small correction to the density matrix:
\begin{align}
i\hbar\partial_{t}\delta\rho=-\frac{\hbar^{2}}{2m}(\Delta_{1}-\Delta_{2})\delta\rho-e(\phi(r_{1},t)-\phi(r_{2},t))\rho_{0}(r_{1}-r_{2}).
\end{align}

Let $$\phi(r,t)=\phi_{\omega k}e^{i(k\cdot r-\omega t)},$$
then the dependence of the solution of (1.1) on the sum $r_{1}+r_{2}$ ( and on the time) can be separated by putting
$$\delta\rho(r_{1},r_{2},t)=\exp[ik\cdot1/2(r_{1}+r_{2})-i\omega t]g_{\omega k}(r_{1}-r_{2}).$$
Substituting this expression in (1.1), we obtain an equation for $g_{\omega k}(R):$
$$\bigg[\hbar\omega+\frac{\hbar^{2}}{2m}(\nabla+\frac{1}{2}ik)^{2}-\frac{\hbar^{2}}{2m}(\nabla-\frac{1}{2}ik)^{2}\bigg]g_{\omega k}(R)
=-e\phi_{\omega k}(e^{ik\cdot R/2}-e^{-ik\cdot R/2})\rho_{0}(R).$$

We can now apply a Fourier expansion with respect to $ R. $ Multiplying both sides by $\exp(-ip\cdot R/\hbar)$ and integrating over $d^{3}x,$ we obtain
$$\bigg[\hbar\omega-\epsilon(p+\frac{1}{2}\hbar k)-\epsilon(p-\frac{1}{2}\hbar k)\bigg]g_{\omega k}(p)
=-(ie\phi_{\omega k}/\hbar k\mathcal{N}_{q})\bigg[n_{0}(p-\frac{1}{2}\hbar k)-n_{0}(p+\frac{1}{2}\hbar k)\bigg],$$
where $\epsilon(p)=p^{2}/2m,$ or equivalently
\begin{align}
g_{\omega k}(p)=(ie\phi_{\omega k}/\hbar k\mathcal{N}_{q})\bigg[n_{0}(p+\frac{1}{2}\hbar k)-n_{0}(p-\frac{1}{2}\hbar k)\bigg]\bigg/(\omega-k\cdot v).
\end{align}

The value of the density matrix at $r_{1}=r_{2}\equiv r$ determines the number density of particles in the system: $N_{q}=2\mathcal{N}_{q}\rho(t,r,r).$ Hence the change in the electron  density due to the field is
$$\delta N_{q}=2\mathcal{N}_{q}\delta\rho(t,r,r)=2\mathcal{N}_{q}e^{i(k\cdot r-\omega t)}g_{\omega k}(R=0),$$
or, expressing $g_{\omega k}(R=0)$ in Fourier component,
\begin{align}
\delta N_{q}=2\mathcal{N}_{q}e^{i(k\cdot r-\omega t)}\int g_{\omega k}(p)d^{3}p/(2\pi\hbar)^{3}.
\end{align}
The corresponding change in the charge density is $-q\delta N_{q}.$

Note that $$-e\delta N_{q}=-\textmd{div} \phi=-ik\cdot \phi,$$
then $$ik\cdot \phi_{\omega k}=2e\mathcal{N}_{q}\int g_{\omega k}(p)d^{3}p/(2\pi\hbar)^{3},$$
and we have \begin{align}
1-\frac{4\pi e^{2}}{\hbar k^{2}}\int\frac{n_{0}(p+\frac{1}{2}\hbar k)-n_{0}(p-\frac{1}{2}\hbar k)}{k\cdot v-\omega-i0}\frac{2d^{3}p}{(2\pi\hbar)^{3}}=0.
\end{align}


Based on the fact that any matter in nature has the wave-particle property,  in this paper we try to study  Landau damping in the quantum physics through the viewpoint and tool used in the classical physics.

 Now we set (analogy to the classical case):
\begin{align}
\epsilon_{l}(\omega,k)-1=-\frac{4\pi e^{2}}{\hbar k^{2}}\int\frac{n_{0}(p+\frac{1}{2}\hbar k)-n_{0}(p-\frac{1}{2}\hbar k)}{k\cdot v-\omega-i0}\frac{2d^{3}p}{(2\pi\hbar)^{3}}.
\end{align}
From the physical viewpoint, $\epsilon_{l}(\omega,k)$ represents the electronic part of the longitudinal permittivity of a plasma having an electron distribution function $n_{0}(p).$ With an appropriate change of integration variable, (1.5) can be written as
$$\epsilon_{l}(\omega,k)=1-\frac{4\pi me^{2}}{k}\int\frac{n_{0}}{(\omega-k\cdot v)^{2}-\hbar^{2}k^{4}/4m^{2}}dv.$$
From (1.5), one can recover the Vlasov-Poisson(VP) dispersion relation by taking the semiclassical limit $\hbar\rightarrow0,$
\begin{align}
\epsilon_{l}(\omega,k)-1=-\frac{4\pi e^{2}}{ k^{2}}\int k\cdot\frac{\partial n_{0}}{\partial p}\frac{1}{k\cdot v-\omega-i0}d^{3}p.
\end{align}

In other words, in the quasi-classical case, the functions $n_{0}(p\pm\frac{1}{2}\hbar k)$ can be expanded in powers of $ k.$ Then
$$n_{0}(p+\frac{1}{2}\hbar k)-n_{0}(p-\frac{1}{2}\hbar k)\approx\hbar k\cdot\partial n_{0}(p)/\partial p,$$
and then (1.5) becomes the quasi-classical case(1.6).

In general, changing the variable of integration in each term in (1.5) by $p\pm\frac{1}{2}\hbar k\rightarrow p,$ we obtain
$$\epsilon_{l}(\omega,k)-1=-\frac{4\pi e^{2}}{\hbar k^{2}}\int\bigg\{\frac{1}{\omega_{+}-k\cdot v+i0}-\frac{1}{\omega_{-}-k\cdot v+i0}\bigg\}n_{0}(p)\frac{2d^{3}p}{(2\pi\hbar)^{3}},$$
where $\omega_{\pm}=\omega\pm\frac{\hbar k^{2}}{2m}.$

We write the complex quantity $\epsilon_{l}=\epsilon_{r}+i\epsilon_{i}.$ Using $\frac{1}{z-i0}=\mathbf{P}\frac{1}{z}+i\pi\delta(z),$ we have
$$\epsilon_{i}=-\frac{4\pi e^{2}}{\hbar k^{2}}\int\bigg\{\delta(\omega_{+}-k\cdot v)-\delta(\omega_{-}-k\cdot v)\bigg\}n_{0}(p)\frac{2d^{3}p}{(2\pi\hbar)^{3}},$$
or $$\epsilon_{i}=-\frac{4m\pi e^{2}}{\hbar k^{2}}\bigg[n_{0}(v)|_{v=\omega_{+}/k}-n_{0}(v)|_{v=\omega_{-}/k}\bigg].$$

From (1.5)-(1.6), it is easy to observe that there are two different points  between the quantum case and the classical case. On the one hand,  the velocity derivative of the distribution function is replaced with a finite difference in the plasmon momentum $\hbar k;$ on the other hand, quantum mechanics can influence the background distribution $n_{0},$ which is a Fermi-Dirac distribution $n_{0}\propto[1+\exp(v^{2}/2v^{2}_{T}-\mu/T)]^{-1}, T$ represents the temperature of the electron.

Since there is no magnetic field associated with the wave, without loss of generality, we can restrict our space direction in the $z-$axis so that the component of $ v $ along it is denoted by $v_{z}.$    Then after integrating over the perpendicular direction,  the corresponding distribution
is \begin{align}
n_{0}(v_{z},\mu)=\frac{N(\mu)}{v_{T}}\ln(1+e^{-v^{2}_{z}/2v^{2}_{T}+\mu/T}),
\end{align}
where $N(\mu)=[\sqrt{\pi}Li_{3/2}(-e^{\mu/T})]^{-1}$ to normalizing $n_{0}$ to 1, and $v_{T}=(T/m)^{1/2}$ is the classical electron thermal velocity. Here, $\mu/T$ determines the level of degeneracy of the system in the quantum state. And if $\mu/T\rightarrow-\infty,$ a Maxwellian in the classical case is recovered.

As in the classical case,  we show that the correct way to perform the integral in the quantum dispersion relation is not simply to take the principal value, but in the complex $v$ plane, following a contour that leaves the singularity always on the same side. With this prescription, the dielectric function is found to possess an  imaginary part, which in turn gives rise to a damping rate $\gamma$ for the wave. And the detailed argument is as follows.

We  assume that $\omega/k\gg v_{T}$ that is the same with that in  the classical case.
Thus, in general, the normal modes will be described by \begin{align}
\epsilon_{l}(\omega,k)=\epsilon_{r}(\omega,k)+i\epsilon_{i}(\omega,k)=0,
\end{align}
where $$\epsilon_{r}(\omega,k)=\epsilon_{r}(\omega_{+},k)-\epsilon_{r}(\omega_{-},k);\quad\quad\epsilon_{i}(\omega,k)=\epsilon_{i}(\omega_{+},k)-\epsilon_{i}(\omega_{-},k).$$

It is usual to assume that
\begin{align}
|\epsilon_{i}|\ll|\epsilon_{r}|,
\end{align}
so that the characteristics of the wave are still primarily determined by $\epsilon_{r}.$ We aim to solve eqn (1.5) for $\omega,$ given $ k. $ It is clear that $ \omega $ will have both a real part and a small imaginary part; let us write
\begin{align}
\omega=\omega_{r}-i\gamma,\quad\quad\quad \omega_{\pm}=\omega_{r}\pm\frac{\hbar k^{2}}{2m}-i\gamma,\quad\quad\quad|\gamma|\ll1,
\end{align}
The time-dependence then becomes
\begin{align}
\exp(-i\omega t)=\exp(-i\omega_{r}t)\exp(-\gamma t),
\end{align}
so that a positive value of $\gamma$ corresponds to damping of a wave that has frequency $\omega_{r}.$ Now consider a Taylor expansion of eqn (1.5), using eqn (1.8) and the inequality eqn(1.9):
\begin{align}
\epsilon_{l}(\omega,k)=\epsilon_{r}(\omega_{r},k)+i\epsilon_{i}(\omega_{r},k)
-i\gamma\bigg(\frac{\partial\epsilon_{r}(\omega,k)}{\partial\omega}\bigg)_{\omega=\omega_{r}}+\ldots=0.
\end{align}

Let us separate real and imaginary parts in eqn (1.12). It follows that the real frequency $\omega_{r}$ is a solution of
\begin{align}
\epsilon_{r}(k,\omega_{r})=0,
\end{align}
and that \begin{align}
\gamma=\frac{\epsilon_{i}(\omega_{r},k)}{\bigg(\frac{\partial\epsilon_{r}(\omega,k)}{\partial\omega}\bigg)_{\omega=\omega_{r}}}.
\end{align}
In fact, eqn (1.14) is a general result, which applies to any relation of the form given by eqn (1.7). For the case in this paper, eqn (1.14) becomes
\begin{align}
\gamma=-\frac{4\pi e^{2}}{\hbar k^{2}\bigg(\frac{\partial\epsilon_{r}(\omega,k)}{\partial\omega}\bigg)_{\omega=\omega_{r}}}\bigg[n_{0}(v_{z})|_{v_{z}=\omega_{+}/k}-n_{0}(v_{z})|_{v_{z}=\omega_{-}/k}\bigg].
\end{align}

then we obtain

$$\epsilon_{i}(k,\omega_{r})=-\frac{4\pi e^{2}}{\hbar k^{2}}\frac{N(\mu)}{v_{T}}\ln\bigg[(1+e^{\frac{\mu}{T}}e^{-\frac{\omega^{2}_{+}}{2k^{2}v_{T}^{2}}})
/(1+e^{\frac{\mu}{T}}e^{-\frac{\omega^{2}_{-}}{2k^{2}v_{T}^{2}}})\bigg]$$

In order to expand the denominator, we state the binomial theorem and two equalities for $kv_{z}/\omega\ll1$:
$$\frac{1}{\omega-kv_{z}}=\frac{1}{\omega}\bigg(1+\frac{kv_{z}}{\omega}+\frac{k^{2}v_{z}^{2}}{\omega^{2}}+\frac{k^{3}v_{z}^{3}}{\omega^{3}}+\cdots\bigg);$$
and for all positive integers $n,$
 $$\int^{\infty}_{-\infty}x^{2n}\exp(-\alpha x^{2})dx=\frac{(2n-1)!}{(2\alpha)^{n}}(\pi/\alpha)^{\frac{1}{2}},$$
 $$\int^{\infty}_{-\infty}x^{2n+1}\exp(-\alpha x^{2})dx=0,$$
 it is easy to check that for any constant $b\neq0,$
 $$\int^{\infty}_{-\infty}(x+b)^{n}\exp(-\alpha x^{2})dx=
 \left\{\begin{array}{l}
 b^{n}\sqrt{\frac{\pi}{\alpha}}+\sum^{\frac{n}{2}}_{k=1}C^{2k}_{n}b^{n-2k}\frac{(2k-1)!}{(2\alpha)^{k}}\sqrt{\frac{\pi}{\alpha}},\quad \textmd{n even}\\
 b^{n}\sqrt{\frac{\pi}{\alpha}}+\sum^{[\frac{n}{2}]}_{k=1}C^{2k}_{n}b^{n-2k}\frac{(2k-1)!}{(2\alpha)^{k}}\sqrt{\frac{\pi}{\alpha}}, \quad \textmd{n odd}.
 \end{array}\right.$$

And  we already know that $$
\epsilon_{r}(k,\omega)=1-\frac{4\pi e^{2}N(\mu)}{\hbar k^{2}v_{T}}\mathbf{P}\int\frac{n_{0}(v_{z}+\frac{1}{2}\hbar k)-n_{0}(v_{z}-\frac{1}{2}\hbar k)}{\omega-kv_{_{z}}}dv_{z},
$$

through simple computation, we have
\begin{align}
&\epsilon_{r}(k,\omega)=1-\frac{4\pi e^{2}N(\mu)e^{\frac{\mu}{T}}}{ kv_{T}\omega}\int(1+kv_{z}/\omega+k^{2}v^{2}_{z}/\omega^{2}+\ldots)\frac{e^{-\frac{(v_{z}-\frac{\hbar k}{2m})^{2}}{2v_{T}^{2}}}-e^{-\frac{(v_{z}+\frac{\hbar k}{2m})^{2}}{2v_{T}^{2}}}}{\hbar k}dv_{z}\notag\\
&=1-\frac{4\pi e^{2}N(\mu)e^{\frac{\mu}{T}}}{ kv_{T}\omega}\int e^{-\frac{v_{z}^{2}}{2v_{T}^{2}}}\frac{(k(v_{z}+\frac{\hbar k}{2m})/\omega+k^{2}(v_{z}+\frac{\hbar k}{2m})^{2}/\omega^{2}+\ldots)-
(k(v_{z}-\frac{\hbar k}{2m})/\omega+k^{2}(v_{z}-\frac{\hbar k}{2m})^{2}/\omega^{2}+\ldots)}{\hbar k}dv_{z}\notag\\
&=1-\frac{4\pi e^{2}N(\mu)e^{\frac{\mu}{T}}}{ kv_{T}\omega}\int e^{-\frac{v_{z}^{2}}{2v_{T}^{2}}}\frac{(k(v_{z}+\frac{\hbar k}{2m})/\omega+k^{2}(v_{z}+\frac{\hbar k}{2m})^{2}/\omega^{2}+\ldots)-
(k(v_{z}-\frac{\hbar k}{2m})/\omega+k^{2}(v_{z}-\frac{\hbar k}{2m})^{2}/\omega^{2}+\ldots)}{\hbar k}dv_{z}\notag\\
&=1-\frac{4\pi e^{2}N(\mu)e^{\frac{\mu}{T}}}{ kv_{T}\omega}\frac{\sum^{\infty}_{n=1}\{[(\frac{(\hbar k)}{2m})^{n}+\sum^{[\frac{n}{2}]}_{l=1}C_{n}^{2l}
(\frac{\hbar k}{2m})^{n-2l}\frac{(2l-1)!}{2^{l/2}}v^{l}_{T}]-[(-\frac{(\hbar k)}{2m})^{n}+\sum^{[\frac{n}{2}]}_{l=1}C_{n}^{2l}
(-\frac{\hbar k}{2m})^{n-2l}\frac{(2l-1)!}{2^{l/2}}v^{l}_{T}]\}}{\hbar k}\notag\\
&=1-\frac{4\pi e^{2}N(\mu)e^{\frac{\mu}{T}}}{ kv_{T}\omega}\bigg(\frac{\frac{\hbar k^{2}v_{T}}{\sqrt{\pi}m\omega}+\frac{\hbar^{3} k^{6}v_{T}}{\pi^{\frac{1}{2}}\omega^{3}(2m)^{3}}+\frac{3\hbar k^{4}v^{2}_{T}}{\pi^{\frac{1}{2}}\omega^{3}\sqrt{2}m}+\ldots}{\hbar k}\bigg)\notag\\
&=1-\frac{\omega^{2}_{pe}}{\omega^{2}}\bigg(1+3\frac{k^{2}v^{2}_{T}}{\omega^{2}}+\frac{\hbar^{2}k^{4}}{4m^{2}\omega^{2}}+\ldots\bigg)\notag\\
&=1-\frac{\omega^{2}_{pe}}{\omega^{2}}\bigg(1+3\frac{k^{2}v^{2}_{T}}{\omega^{2}}+\frac{\hbar^{4}k^{4}}{4m^{2}\hbar^{2}\omega^{2}}+\ldots\bigg)
\end{align}
where $\omega_{pe}=(n_{0}e^{2}/m\epsilon_{0})^{\frac{1}{2}}$ is the plasma frequency.

In this paper we restrict our attention to the case where the value of $v_{z}$ that satisfies $v_{z}=\omega/k$ exceeds $v_{T},$ so that the number of the resonant electrons is exponentially small and their effect is correspondingly weak. Then the range of values of $v_{z}$ of interest is $kv_{z}/\omega\ll1.$

Then we know that to leading order $\omega_{\pm r}=\omega_{pe}\pm\frac{\hbar k^{2}}{2m},$ so that differentiating eqn (1.13) with respect to $\omega$ yields, through good approximation:

  \begin{align}
 \bigg( \frac{\partial\epsilon_{r}(\omega,k)}{\partial\omega}\bigg)_{\omega=\omega_{r}}=\frac{2}{\omega_{pe}},
  \end{align}
and \begin{align}
\gamma=-\frac{8\pi e^{2}}{\hbar k^{2}\omega_{pe}}\frac{N(\mu)}{v_{T}}\ln\bigg[(1+e^{\frac{\mu}{T}}e^{-\frac{\omega^{2}_{+}}{2k^{2}v_{T}^{2}}})
/(1+e^{\frac{\mu}{T}}e^{-\frac{\omega^{2}_{-}}{2k^{2}v_{T}^{2}}})\bigg].
\end{align}

If  assuming $-\mu/T$ large enough, then \begin{align}
\gamma\approx-\frac{8\pi e^{2}}{\hbar k^{2}v_{T}\omega_{pe}}e^{-\frac{\omega^{2}_{pe}}{2k^{2}v^{2}_{T}}-\frac{\hbar^{2}k^{2}}{4m^{2}v_{T}^{2}}}
\bigg(e^{-\frac{\hbar\omega_{pe}}{mv^{2}_{T}}}-e^{\frac{\hbar\omega_{pe}}{mv^{2}_{T}}}\bigg).
\end{align}

At the same time, we assume that $\frac{\hbar k}{2m}\ll v_{T},$ then $\gamma\approx e^{-\frac{\omega^{2}_{pe}}{2k^{2}v^{2}_{T}}},$
that is the same to  Landau damping rate in classical plasma.

           \begin{center}
\item\section{ The mathematical tool}
\end{center}

        Now we introduce some notations. We denote $\mathbb{T}^{3}=\mathbb{R}^{3}/\mathbb{Z}^{3}.$ For function $f(x,v),$ we define the Fourier transform as follows.

        For a function $f=f(x),x\in\mathbb{T}^{d},$ we define its Fourier transform as follows:
 $$\hat{f}(k)=\frac{1}{(2\pi)^{d}}\int_{\mathbb{T}^{d}}f(x)e^{-ix\cdot k}dx,\quad k\in\mathbb{Z}^{d}.$$
 Similarly, for a function $f=f(v),v\in\mathbb{R}^{d},$ we define its Fourier transform by:
 $$\hat{f}(\xi)=\frac{1}{(2\pi)^{d}}\int_{\mathbb{R}^{d}}f(v)e^{-iv\cdot \xi}dv,\quad \xi\in\mathbb{R}^{d}.$$

 Finally, if $f=f(x,v), (x,v)\in\mathbb{T}^{d}\times\mathbb{R}^{d},$ we define its Fourier transform through the following formula:
 $$\hat{f}(k,\xi)=\frac{1}{(2\pi)^{d}}\int_{\mathbb{T}^{d}\times\mathbb{R}^{d}}f(x,v)e^{-ik\cdot x-iv\cdot\xi}dxdv,\quad (k,\xi)\in\mathbb{Z}^{d}\times\mathbb{R}^{d}.$$

 We shall also use the Fourier transform in time, if $f=f(t),t\in\mathbb{R},$ we denote
 $$\tilde{f}(\omega)=\int_{\mathbb{R}}f(t)e^{-it\omega}dt,\quad \omega\in\mathbb{C}.$$

Now we start to introduce the very important tools in our paper. These are time-shift pure and hybrid analytic norms. They are the same with those
           in the paper [28] written by Mouhot and Villani.
            \begin{de}(Hybrid analytic norms)
            $$\|f\|_{\mathcal{C}^{\lambda,\mu}}=\sum_{m,n\in\mathbb{N}_{0}^{3}}\frac{\lambda^{n}}{n!}\frac{\mu^{m}}{m!}
            \|\nabla^{m}_{x}\nabla_{v}^{n}f\|_{L^{\infty}(\mathbb{T}_{x}^{3}\times\mathbb{R}^{3}_{v})},\quad
            \|f\|_{\mathcal{F}^{\lambda,\mu}}=\sum_{k\in\mathbb{Z}^{3}}\int_{\mathbb{R}^{3}}|\tilde{f}(k,\eta)|
           e^{2\pi\lambda|\eta|}e^{2\pi\mu|k|}d\eta,$$
            $$\|f\|_{\mathcal{Z}^{\lambda,\mu}}=\sum_{l\in\mathbb{Z}^{3}}\sum_{n\in\mathbb{N}_{0}^{3}}
            \frac{\lambda^{n}}{n!}e^{2\pi\mu|l|}\|\widehat{\nabla_{v}^{n}f(l,v)}\|_{L^{\infty}(\mathbb{R}_{v}^{3})}.$$
            \end{de}
           \begin{de} (Time-shift pure and hybrid analytic norms)
For any $\lambda,\mu\geq0,p\in[1,\infty],$ we define
            $$\|f\|_{\mathcal{C}^{\lambda,\mu}_{\tau}}=\|f\circ S^{0}_{\tau}(x,v)\|_{\mathcal{C}^{\lambda,\mu}}=\sum_{m,n\in\mathbb{N}_{0}^{3}}\frac{\lambda^{n}}{n!}\frac{\mu^{m}}{m!}
            \|\nabla^{m}_{x}(\nabla_{v}+\tau\nabla_{x})^{n}f\|_{L^{\infty}(\mathbb{T}_{x}^{3}\times\mathbb{R}^{3}_{v})},$$
           $$\|f\|_{\mathcal{F}^{\lambda,\mu}_{\tau}}=\|f\circ S^{0}_{\tau}(x,v)\|_{\mathcal{F}^{\lambda,\mu}}=\sum_{k\in\mathbb{Z}^{3}}\int_{\mathbb{R}^{3}}|\tilde{f}(k,\eta)|
           e^{2\pi\lambda|k\tau+\eta|}e^{2\pi\mu|k|}d\eta,$$
         $$\|f\|_{\mathcal{Z}^{\lambda,\mu}_{\tau}}=\|f\circ S^{0}_{\tau}(x,v)\|_{\mathcal{Z}^{\lambda,\mu}}=\sum_{l\in\mathbb{Z}^{3}}\sum_{n\in\mathbb{N}_{0}^{3}}
            \frac{\lambda^{n}}{n!}e^{2\pi\mu|l|}\|(\nabla_{v}+2i\pi\tau\cdot l)^{n}\hat{f}(l,v)\|
            _{L^{\infty}(\mathbb{R}_{v}^{3})},$$
               $$\|f\|_{\mathcal{Z}^{\lambda,\mu;p}_{\tau}}=\sum_{l\in\mathbb{Z}^{3}}\sum_{n\in\mathbb{N}_{0}^{3}}
            \frac{\lambda^{n}}{n!}e^{2\pi\mu|l|}\|(\nabla_{v}+2i\pi\tau\cdot l)^{n}\hat{f}(l,v)\|_{L^{p}(\mathbb{R}_{v}^{3})},$$
            $$\|f\|_{\mathcal{Y}^{\lambda,\mu}_{\tau}}=\|f\|_{\mathcal{F}^{\lambda,\mu;\infty}_{\tau}}
            =\sup_{k\in\mathbb{Z}^{3},\eta\in\mathbb{R}^{3}}e^{2\pi\mu|k|}e^{2\pi\lambda|\eta+k\tau|}|\hat{f}(k,\eta)|.$$
           \end{de}

           From the above  definitions, we can state some simple and important propositions, and the related proofs can be found in [26.28], so we remove the proofs.
\begin{prop} For any $\tau\in\mathbb{R},\lambda,\mu\geq 0,$
\begin{itemize}
 \item[{(i)}]   if $f$ is a function only of $x,$ then
           $\|f\|_{\mathcal{C}^{\lambda,\mu}_{\tau}}=\|f\|_{\mathcal{C}^{\lambda|\tau|+\mu}},\|f\|_{\mathcal{F}^{\lambda,\mu}_{\tau}}
           =\|f\|_{\mathcal{Z}^{\lambda,\mu}_{\tau}}=\|f\|_{\mathcal{F}^{\lambda|\tau|+\mu}};$
          \item[{(ii)}]    if $f$ is a function only of $v,$ then
           $\|f\|_{\mathcal{C}^{\lambda,\mu;p}_{\tau}}=\|f\|_{\mathcal{Z}^{\lambda,\mu;p}_{\tau}}=\|f\|_{\mathcal{C}^{\lambda,;p}};$
            \item[{(iii)}] for any $\lambda>0,$ then
            $\|f\circ(Id+G)\|_{\mathcal{F}^{\lambda}}\leq\|f\|_{\mathcal{F}^{\lambda+\nu}},\nu=\|G\|_{\dot{\mathcal{F}}^{\lambda}};$
            \item[{(iv)}] for any $\bar{\lambda}>\lambda,p\in[1,\infty],$ $\|\nabla f\|_{\mathcal{C}^{\lambda;p}}\leq
            \frac{1}{\lambda e\log(\frac{\bar{\lambda}}{\lambda} )}
            \| f\|_{\mathcal{C}^{\bar{\lambda};p}},$
            $\|\nabla f\|_{\mathcal{F}^{\lambda;p}}\leq\frac{1}{2\pi e(\bar{\lambda}-\lambda )}
            \| f\|_{\mathcal{F}^{\bar{\lambda};p}},$
            \item[{(v)}]  for any $\bar{\lambda}>\lambda>0,\mu>0,$ then
            $\|vf\|_{\mathcal{Z}^{\lambda,\mu;1}_{\tau}}\leq\|f\|_{\mathcal{Z}^{\bar{\lambda},\mu;1}_{\tau}};$
           \item[{(vi)}]  for any $\bar{\lambda}>\lambda,\bar{\mu}>\mu,$
          $\|\nabla_{v}f\|_{\mathcal{Z}_{\tau}^{\lambda,\mu;p}}
          \leq C\bigg(\frac{1}{\lambda\log(\frac{\bar{\lambda}}{\lambda})}\|f\|_{\mathcal{Z}_{\tau}^{\bar{\lambda},\bar{\mu};p}}+\frac{\tau}{\bar{\mu}-\mu}
 \|f\|_{\dot{\mathcal{Z}}_{\tau}^{\bar{\lambda},\bar{\mu};p}}\bigg);$
 \item[{(vii)}] for any $\bar{\lambda}>\lambda,$
 $\|(\nabla_{v}+\tau\nabla_{x})f\|_{\mathcal{Z}_{\tau}^{\lambda,\mu;p}}\leq \frac{1}{C\lambda\log(\frac{\bar{\lambda}}{\lambda})}
 \|f\|_{\mathcal{Z}_{\tau}^{\bar{\lambda},\mu;p}};$
  \item[{(viii)}] for any $\bar{\lambda}\geq\lambda\geq0,\bar{\mu}\geq\mu\geq0,$ then $\|f\|_{\mathcal{Z}_{\tau}^{\lambda,\mu}}
  \leq_{\mathcal{Z}_{\tau}^{\bar{\lambda},\bar{\mu}}}.$ Moreover, for any $\tau,\bar{\tau}\in\mathbb{R},$ $p\in[1,\infty],$ we have
  $\|f\|_{\mathcal{Z}_{\tau}^{\lambda,\mu;p}}\leq\|f\|_{\mathcal{Z}_{\bar{\tau}}^{\lambda,\mu+\lambda|\tau-\bar{\tau}|;p}};$
  \item[{(viiii)}] $\|f\|_{\mathcal{Y}^{\lambda,\mu}_{\tau}}\leq\|f\|_{\mathcal{Z}^{\lambda,\mu;1}_{\tau}};$
  \item[{(ix)}] for any function $f=f(x,v),$ $\|\int_{\mathbb{R}^{3}}fdv\|_{\mathcal{F}^{\lambda|\tau|+\mu}}\leq\|f\|_{\mathcal{Z}^{\lambda,\mu;1}_{\tau}}.$
 \end{itemize}
\end{prop}

           \begin{prop} For any $X\in\{\mathcal{C},\mathcal{F},\mathcal{Z}\}$ and any $t,\tau\in\mathbb{R},$
           $$\|f\circ S^{0}_{\tau}\|_{X^{\lambda,\mu}_{\tau}}=\|f\|_{X^{\lambda,\mu}_{t+\tau}}.$$

           \end{prop}
           \begin{lem} Let $\lambda,\mu\geq0,t\in\mathbb{R},$ and consider two functions $F,G:\mathbb{T}^{3}\times\mathbb{R}^{3}\rightarrow\mathbb{T}^{3}\times\mathbb{R}^{3}.$
           Then there is $ \varepsilon\in(0,\frac{1}{2})$ such that  if $F,G$ satisfy
\begin{align}
\|\nabla(F-Id)\|_{\mathcal{Z}^{\lambda,\mu}_{\tau}}\leq\varepsilon,
\end{align}
where $\lambda=\lambda +2\|F-G\|_{\mathcal{Z}^{\lambda ,\mu}_{\tau}},\quad\mu=\mu+2(1+|\tau|)\|F-G\|_{\mathcal{Z}^{\lambda ,\mu}_{\tau}} ,$ then $F$ is invertible  and
\begin{align}
&\|F^{-1}\circ G-Id\|_{\mathcal{Z}^{\lambda ,\mu}_{\tau}}\leq2\|F-G\|_{\mathcal{Z}^{\lambda ,\mu}_{\tau}}.
\end{align}
\end{lem}

\begin{prop} For any $\lambda,\mu\geq0$ and any $p\in[1,\infty],\tau\in\mathbb{R},\sigma\in\mathbb{R},a\in\mathbb{R}\setminus\{0\}$ and $b\in\mathbb{R},$
we have
$$\|f(x+bv+X(x,v),av+V(x,v))\|_{\mathcal{Z}^{\lambda,\mu;p}_{\tau}}\leq|a|^{-\frac{3}{p}}\|f\|_{\mathcal{Z}^{\alpha,\beta;p}_{\sigma}},$$
where $\alpha=\lambda|a|+\|V\|_{\mathcal{Z}^{\lambda,\mu}_{\tau}},\quad \beta=\mu+\lambda|b+\tau-a\sigma|+\|X-\sigma V\|_{\mathcal{Z}^{\lambda,\mu}_{\tau}}.$
\end{prop}
\begin{lem} Let $G=G(x,v)$ and $R=R(x,v)$ be valued in $\mathbb{R},$ and
$\beta(x)=\int_{\mathbb{R}^{3}}(G\cdot R)(x,v)dv.$
Then for any $\lambda,\mu,t\geq0$ and any $b>-1,$ we have
$$\|\beta\|_{\mathcal{F}^{\lambda t+\mu}}\leq3\|G\|_{\mathcal{Z}_{t-\frac{bt}{1+b}}^{\lambda(1+b),\mu;1}}
\|R\|_{\mathcal{Z}_{t-\frac{bt}{1+b}}^{\lambda(1+b),\mu}}.$$
\end{lem}
$$$$
\begin{center}
\item\section{The main mathematical  properties }
\end{center}

\begin{center}
\item\subsection{ Mathematical expression of the Kernel }
\end{center}

From the above physical point of view, under the assumption of the stability condition, we have known that,  scattering resonances occurring
 at distinct frequencies are asymptotically well separated. In the following, through  complicated computation, we give a detailed description by
  using mathematical tool. And the proof is simple, so we omit.

 \begin{thm} Let $\lambda,\bar{\lambda},\mu,\bar{\mu},\mu'$ be such that
$2\lambda\geq\bar{\lambda}>\lambda>0, \bar{\mu}\geq\mu'>\mu>0,$ and let $\gamma\geq0$ and $b=b(t,s)>0,$
  $ R=R(t,x),G=G(t,x,v)$ and assume $\widehat{G}(t,0,v)=0,$ we have, if
$$\sigma(t,x,v)=\int^{t}_{0}R(s,x+(t-s)v)G(s,x+(t-s)v,v)ds,$$
 then
\begin{align}
 & \|\sigma(t,\cdot)\|_{\mathcal{Z}^{\lambda,\mu;1}_{t}}\leq\int^{t}_{0} K(t,s)\|R\|_{\mathcal{F}^{\lambda s+\mu'-\lambda b(t-s)}}\frac{\|G\|_{\mathcal{Z}^{\bar{\lambda}(1+b),\bar{\mu};1}_{s-\frac{bt}{1+b}}}}{1+s}ds,
\end{align}
where $K(t,s)=(1+s)\sup_{k,l\in\mathbb{Z}^{3}_{\ast}}\frac{e^{-\pi(\bar{\mu}-\mu)|l|}
e^{-\pi(\bar{\lambda}-\lambda)|k(t-s)+ls|}e^{-2\pi[\mu'-\mu+\lambda b(t-s)]|k-l|}}{1+|k-l|^{\gamma}}.$

\end{thm}

\begin{center}
\item\subsection{The  time evolution of the Kernel}
\end{center}

 In the following we  have to analyze  the influence of scattering resonance when studying the stability of  particles as $t\rightarrow\infty.$ This is similar to the occurrence  of Landau damping in the classical plasma and the Diophantus condition of KAM theory.
\begin{prop}(Exponential moments of the kernel) Let $\gamma\in[1,\infty)$ be given. For any $\alpha\in(0,1),$ let $K^{(\alpha),\gamma}$ be defined
$$K^{(\alpha),\gamma}(t,s)=(1+s)\sup_{k,l\in\mathbb{Z}^{3}_{\ast}}\frac{e^{-\alpha|l|}e^{-\alpha(t-s)\frac{|k-l|}{t}}e^{-\alpha|k(t-s)+ls|}}{1+|k-l|^{\gamma}}.$$
Then for any $\gamma<\infty,$ there is $\bar{\alpha}=\bar{\alpha}(\gamma)>0$ such that if $\alpha\leq\bar{\alpha}$ and $\varepsilon\in(0,1),$ then for any
$t>0,$
$$e^{-\varepsilon t}\int^{t}_{0}K^{(\alpha),\gamma}(t,s)e^{\varepsilon s}ds
\leq C\bigg(\frac{1}{\alpha\varepsilon^{\gamma}t^{\gamma-1}}+\frac{1}{\alpha\varepsilon^{\gamma}t^{\gamma}}\log\frac{1}{\alpha}
+\frac{1}{\alpha^{2}\varepsilon^{1+\gamma}t^{1+\gamma}}+\bigg(\frac{1}{\alpha^{3}}+\frac{1}{\alpha^{2}\varepsilon}\log\frac{1}{\alpha}\bigg)
e^{-\frac{\varepsilon t}{4}}+\frac{e^{-\frac{\alpha t}{2}}}{\alpha^{3}}\bigg),$$
where $C=C(\gamma).$

In particular, if $\gamma>1,\varepsilon\leq\alpha,$ then
$e^{-\varepsilon t}\int^{t}_{0}K^{(\alpha),\gamma}(t,s)e^{\varepsilon s}ds\leq\frac{C(\gamma)}{\alpha^{3}\varepsilon^{1+\gamma}t^{\gamma-1}};$

if $\gamma=1,\varepsilon\leq\alpha,$ then
$e^{-\varepsilon t}\int^{t}_{0}K^{(\alpha),\gamma}(t,s)e^{\varepsilon s}ds\leq\frac{C}{\alpha^{3}}(\frac{1}{\varepsilon}+\frac{1}{\varepsilon^{2}t}).$
\end{prop}

\begin{prop} With the same notation as in Proposition 3.2, for any $\gamma\geq1,$ we have
\begin{align}
\sup_{s\geq0}e^{\varepsilon s}\int^{\infty}_{s}e^{-\varepsilon t}K^{(\alpha),\gamma}(t,s)dt\leq C(\gamma)\bigg(\frac{1}{\alpha^{2}\varepsilon}+
\frac{1}{\alpha\varepsilon^{\gamma}}\log\frac{1}{\alpha}\bigg)
\end{align}
\end{prop}


\begin{center}
\item\subsection{ The ``Maximal" principle}
\end{center}

From now on, we will state the main result of this section that is the same with section 7.4 in [26,28], the detailed proof can be found in appendix ( also see [26,28] ).
We define $\|\Phi(t)\|_{\lambda}=\sum_{k\in\mathbb{Z}_{\ast}^{3}}|\Phi(k,t)|e^{2\pi\lambda|k|}.$

\begin{thm} Assume that $w_{\hbar}^{0}(v),W=W(x)$ satisfy the conditions of Theorem 0.1, and the Stability condition holds.
Let $A\geq0,\mu\geq0$ and $\lambda\in(0,\lambda^{\ast}]$
with $0<\lambda^{\ast}<\lambda_{0}.$ Let $(\Phi(k,t))_{k\in\mathbb{Z}_{\ast}^{3},t\geq0}$ be a continuous functions of $t\geq0,$ valued in
$\mathbb{C}^{\mathbb{Z}^{3}_{\ast}},$ such that for all $t\geq0,$
\begin{align}\|\Phi(t)-\int^{t}_{0}K^{0}(t-s)\Phi(s)ds\|_{\lambda t+\mu}\leq A+\int^{t}_{0}(K_{0}(t,s)+K_{1}(t,s)+\frac{c_{0}}{(1+s)^{m}})
\|\Phi(s)\|_{\lambda s+\mu} ds,
\end{align}
where $c_{0}\geq0,m>1,$ and $K_{0}(t,s),K_{1}(t,s)$ are non-negative kernels. Let $\varphi(t)=\|\Phi(t)\|_{\lambda t+\mu}.$
Then we have the following:

(i) Assume that $\gamma>1$ and $K_{1}=cK^{(\alpha),\gamma}$ for some $c>0,\alpha\in(0,\bar{\alpha}(\gamma)),$ where $K^{(\alpha),\gamma},\bar{\alpha}(\gamma)$
are the same with that
 defined by Proposition 3.2. Then there are positive constants $C,\chi,$ depending only on $\gamma,\lambda^{\ast},\lambda_{0},\kappa,c_{0}, C_{W}$ and $ m,$ uniform
 as $\gamma\rightarrow1,$ such that if $\sup_{t\geq0}\int^{t}_{0}K_{0}(t,s)ds\leq\chi$ and $\sup_{t\geq0}$
 $\bigg(\int^{t}_{0}K_{0}(t,s)^{2}ds\bigg)^{\frac{1}{2}}+\sup_{t\geq0}\int^{\infty}_{t}K_{0}(t,s)dt\leq1,$ then for any $\varepsilon\in(0,\alpha),$ for all $t\geq0,$
 \begin{align}\varphi(t)\leq CA\frac{1+c^{2}_{0}}{\sqrt{\varepsilon}}e^{Cc_{0}}\bigg(1+\frac{c}{\alpha\nu}\bigg)e^{CT}e^{Cc(1+T^{2})}e^{\varepsilon t}
 \end{align}
 where $T_{\nu}=C\max\bigg\{\bigg(\frac{c^{2}}{\alpha^{5}}\varepsilon^{2+\gamma}\bigg)^{\frac{1}{\gamma-1}},
\bigg(\frac{c}{\alpha^{2}}\varepsilon^{\frac{1}{2}+\gamma}\bigg)^{\frac{1}{\gamma-1}},\bigg(\frac{c^{2}_{0}}{\varepsilon}\bigg)^{\frac{1}{2m-1}}\bigg\}.$

(ii) Assume that $K_{1}=\sum^{N}_{j=1}c_{j}K^{(\alpha_{j},1)}$ for some $\alpha_{j}\in(0,\bar{\alpha}(\gamma)),$ where $\bar{\alpha}(\gamma)$ also appears in proposition 3.2;
then there is a numeric constant $\Gamma>0$ such that whenever $1\geq\varepsilon\geq\Gamma\sum^{N}_{j=1}\frac{c_{j}}{\alpha^{3}_{j}},$ with the same notation
 as in (i), for all $t\geq0,$  one has,

  \begin{align}\varphi(t)\leq CA\frac{1+c^{2}_{0}}{\sqrt{\varepsilon}}e^{Cc_{0}}e^{CT}e^{Cc(1+T^{2})}e^{\varepsilon t}
 \end{align}
 where $c=\sum^{N}_{j=1}c_{j}$ and $T=\max\bigg\{\frac{1}{\varepsilon^{2}}\sum^{N}_{j=1}\frac{c_{j}}{\alpha^{3}_{j}},
 \bigg(\frac{c^{2}_{0}}{\varepsilon}\bigg)^{\frac{1}{2m-1}}\bigg\}.$

\end{thm}

$$$$
\begin{center}
\item\section{The property of  the Wigner function}
\end{center}

\begin{center}
\item\subsection{The property of  the Wigner function in the linear case}
\end{center}

In this section we first consider the Wigner function $w_{\hbar}$ defined in (0.11) where $\psi(t,x)$ satisfies (0.1)-(0.3), satisfying the following equations:
\begin{align}
\left\{\begin{array}{l}\partial_{t}w_{\hbar}+v\cdot\nabla_{x}w_{\hbar}+\frac{e}{m}\Theta_{\hbar}[\phi]w^{0}_{\hbar}=0,,\\
(\Theta_{\hbar}[\phi]w^{0}_{\hbar})(t,x,v)
=\frac{im}{(2\pi)^{3}}\int_{\mathbb{R}^{3}}\int_{\mathbb{R}^{3}}\frac{\phi(t,x+\frac{\hbar}{2m}\eta)-\phi(t,x-\frac{\hbar}{2m}\eta)}{\hbar}w^{0}_{\hbar}(t,x,v')e^{i(v-v')\cdot\eta}
dv'd\eta,\\
\phi(t,x)=W(x)\ast n_{\hbar}(t,x),\\
W(x)=W(-x),\widehat{W}(k)\leq\frac{C_{W}}{|k|^{1+\gamma}},\\
w_{\hbar}(0,x,v)=w_{\hbar I}(x,v).
\end{array}\right.
\end{align}

\begin{thm}
Consider equations (4.1). For any $\eta,v\in\mathbb{R}^{3},k\in\mathbb{N}_{0}^{3},$  assume that the following conditions hold for $\bar{\lambda}>\lambda>0, \bar{\mu}>\mu'>\mu>0$:
\begin{itemize}
        \item[(i)] $ \widehat{W}(0)=0$ where $|\widehat{W}(k)|\leq\frac{1}{1+|k|^{\gamma}},\gamma\geq1;$
         \item[(ii)] $||w_{\hbar}^{0}||_{\mathcal{C}^{\lambda;1}}\leq C_{0},$
              for some constants $\lambda,C_{0}>0;$

         \item[(iii)]$||w_{\hbar I}||_{\mathcal{Z}^{\lambda,\mu;1}}\leq \delta_{0}$ for some constants $\mu>0,\delta_{0}>0$;

             \item[(iv)] $\mathbf{ Stability}$ condition: there is some small enough constant $\delta_{0}>0,$  \begin{align}
\bigg\|\int_{\mathbb{T}^{3}}\delta w^{0}_{\hbar}(x,\cdot)dx\bigg\|_{\mathcal{C}^{\bar{\lambda}(1+b);1}}\leq \delta_{0},\quad\quad
\sup_{0<\tau\leq t}\|\delta w^{0}_{\hbar}\|_{\mathcal{Z}^{\bar{\lambda}(1+b),\bar{\mu};1}_{\tau-\frac{bt}{1+b}}}\leq\delta_{0},
\end{align}
         \end{itemize}
         then 
for any fixed $\eta,k,$  we have
         \begin{align}
         &|\hat{w}_{\hbar}(t,k,\eta)-\hat{w}_{\hbar I}(k,\eta)|\leq C(C_{0},\delta_{0}) e^{-2\pi\lambda|\eta+kt|}e^{-2\pi\mu |k|},\notag\\
         &|\hat{n}_{\hbar}(t,k)-\hat{n}_{\hbar I}|\leq C(C_{0},\delta_{0}) e^{-2\pi(\lambda t+\mu)|k|},\notag\\
        & |\hat{\phi}(t,k)|\leq C(C_{0},\delta_{0})e^{-2\pi\lambda|k|t}e^{-2\pi\mu|k| },
         \end{align}
        where $n_{\hbar I}=\int_{\mathbb{T}^{3}}\int_{\mathbb{R}^{3}}w_{\hbar I}(x,v)dvdx.$
         \end{thm}

$ \mathbf{proof}$: Using the Duhamle principle,
\begin{align}
w_{\hbar}(t,x,v)=w_{\hbar I}(x-vt,v)-\int^{t}_{0}[(\Theta_{\hbar}[\phi]w^{0}_{\hbar})](\tau,x-v(t-\tau),v)d\tau,
\end{align}
and \begin{align}
n_{\hbar}(t,x)=\int_{\mathbb{R}^{3}}w_{\hbar I}(x-vt,v)dv-\int^{t}_{0}\int_{\mathbb{R}^{3}}[\Theta_{\hbar}[\phi]w^{0}_{\hbar}](\tau,x-v(t-\tau),v)dvd\tau.
    \end{align}
    Hence, for $k\neq0, k\neq l,$ we have
   $$ \hat{n}_{\hbar}(t,k)=\int_{\mathbb{T}^{3}}\int_{\mathbb{R}^{3}}w_{\hbar I}(x-vt,v)e^{ix\cdot k}dvdx$$
   $$ -\frac{im}{(2\pi)^{3}}\int^{t}_{0}\int_{\mathbb{T}^{3}}\int_{\mathbb{R}^{3}}\int_{\mathbb{R}^{3}}\int_{\mathbb{R}^{3}}
    \frac{\phi(t,x+\frac{\hbar}{2m}\eta)-\phi(t,x-\frac{\hbar}{2m}\eta)}{\hbar}w^{0}_{\hbar}(x,v')e^{i(v-v')\cdot\eta}e^{iv(t-\tau)\cdot k}e^{ix\cdot k}dv'd\eta dxdvd\tau$$
    $$=\widehat{w}_{\hbar I}(k,kt) -\frac{im}{(2\pi)^{3}}\int^{t}_{0}\int_{\mathbb{T}^{3}}\int_{\mathbb{R}^{3}}\int_{\mathbb{R}^{3}}
    \frac{\phi(t,x+\frac{\hbar}{2m}\eta)-\phi(t,x-\frac{\hbar}{2m}\eta)}{\hbar}w^{0}_{\hbar}(x,\eta)e^{iv\cdot\eta}e^{iv(t-\tau)\cdot k}e^{ix\cdot k}d\eta dxdvd\tau$$
     $$=\widehat{w}_{\hbar I}(k,kt) -\frac{im}{(2\pi)^{3}}\int^{t}_{0}\int_{\mathbb{R}^{3}}(k-l)\widehat{\phi}(t,k-l)
    \frac{w^{0}_{\hbar}(l,v-\frac{\hbar}{2m}(k-l))-w^{0}_{\hbar}(l,v+\frac{\hbar}{2m}(k-l))}{\hbar (k-l)}e^{iv(t-\tau)\cdot k}dvd\tau\footnote{The term in the integral term depending on $\hbar$ implies that the quantum damping rate depends on $\hbar$.}$$
    Claim 1: \begin{align}
    \int_{\mathbb{T}^{3}}n_{\hbar}(t,x)dx=\int_{\mathbb{T}^{3}}\int_{\mathbb{R}^{3}}w_{\hbar I}(x-vt,v)dvdx
    \end{align}

    So we have $$\|n_{\hbar}(t,\cdot)\|_{\mathcal{F}^{\lambda t+\mu}}$$
    \begin{align}
   \leq\bigg\|\int_{\mathbb{R}^{3}}w_{\hbar I}(x-vt,v)dv\bigg\|_{\mathcal{F}^{\lambda t+\mu}}
   +\int^{t}_{0} K_{1}(t,\tau)\|n_{\hbar}\|_{\mathcal{F}^{\lambda \tau+\mu'-\lambda b(t-\tau);\gamma}}\frac{\|\delta w^{0}_{\hbar}\|_{\mathcal{Z}^{\bar{\lambda}(1+b),\bar{\mu};1}_{\tau-\frac{bt}{1+b}}}}{1+\tau}d\tau
\end{align}
$$+3\int^{t}_{0}e^{-\pi(\bar{\lambda}-\lambda)(t-\tau)}\bigg\|\int_{\mathbb{T}^{3}}\delta w^{0}_{\hbar}(x,\cdot)dx\bigg\|_{\mathcal{C}^{\bar{\lambda}(1+b);1}}
\|n_{\hbar}\|_{\mathcal{F}^{\lambda \tau+\mu'-\lambda b(t-\tau)}}d\tau$$
where $$K_{1}(t,\tau)=(1+\tau)\sup_{k,l\in\mathbb{Z}^{3}_{\ast}}e^{-\pi(\bar{\mu}-\mu)|l|}
e^{-\pi(\bar{\lambda}-\lambda)|k(t-\tau)+l\tau|}e^{-2\pi[\mu'-\mu+\lambda b(t-\tau)]|k-l|}.$$

By the Stability condition, there is some small enough constant $\delta_{0}>0,$  \begin{align}
\bigg\|\int_{\mathbb{T}^{3}}\delta w^{0}_{\hbar}(x,\cdot)dx\bigg\|_{\mathcal{C}^{\bar{\lambda}(1+b);1}}\leq \delta_{0},\quad\quad
\|\delta w^{0}_{\hbar}\|_{\mathcal{Z}^{\bar{\lambda}(1+b),\bar{\mu};1}_{\tau-\frac{bt}{1+b}}}\leq\delta_{0},
\end{align}
    then we can get
    $$\|n_{\hbar}(t,\cdot)\|_{\mathcal{F}^{\lambda t+\mu}}$$
    \begin{align}
    \leq\bigg\|\int_{\mathbb{R}^{3}}w_{\hbar I}(x-vt,v)dv\bigg\|_{\mathcal{F}^{\lambda t+\mu}} +\int^{t}_{0} K_{1}(t,\tau)\|n_{\hbar}(t,\cdot)\|_{\mathcal{F}^{\lambda \tau+\mu'-\lambda b(t-\tau)}}\frac{\delta_{0}}{1+\tau}d\tau
    \end{align}
$$+\int^{t}_{0}e^{-\pi(\bar{\lambda}-\lambda)(t-\tau)}\delta_{0}
\|n_{\hbar}(t,\cdot)\|_{\mathcal{F}^{\lambda \tau+\mu'-\lambda b(t-\tau)}}d\tau$$
and applying Propositions 3.1-3.4, \begin{align}
\|n_{\hbar}(t,\cdot)\|_{\mathcal{F}^{\lambda t+\mu}}\leq\frac{C\varepsilon_{0}}{\sqrt{\varepsilon}(\lambda-\bar{\lambda})^{2}}
\bigg(1+\frac{1}{\alpha \varepsilon}\bigg)e^{CT}e^{\varepsilon t}e^{Cc(1+T^{2})}
\end{align}
where $\alpha=\bar{\lambda}-\lambda,T=C\max\bigg\{\bigg(\frac{c^{2}}{\alpha^{5}\varepsilon^{2+\gamma}}\bigg)^{1/\gamma-1},
\bigg(\frac{c}{\alpha^{2}\varepsilon^{\gamma+\frac{1}{2}}}\bigg)^{1/\gamma-1}\bigg\}.$

At the same time,
 $$\|w_{\hbar}(t,\cdot)\|_{\mathcal{Z}^{\lambda',\mu';1}_{t}}\leq\|w_{\hbar I}(x-vt,v)\|_{\mathcal{Z}^{\lambda',\mu';1}_{t}}+\int^{t}_{0}\|w^{0}_{\hbar }(x-v\tau,v)\|_{\mathcal{Z}^{\lambda',\mu';1}_{\tau}}\|\Theta_{\hbar}[\phi]\|_{\mathcal{Z}^{\lambda',\mu'}_{\tau}}d\tau$$
 \begin{align}
 \leq\|w_{\hbar I}(x-vt,v)\|_{\mathcal{Z}^{\lambda',\mu';1}_{t}}+\int^{t}_{0}\|w^{0}_{\hbar }(x-v\tau,v)\|_{\mathcal{Z}^{\lambda',\mu';1}_{\tau}}\|n_{\hbar}(\tau,\cdot)\|_{\mathcal{Z}^{\lambda',\mu'}_{\tau}}d\tau
 \end{align}
 $$=\|w_{\hbar I}(x-vt,v)\|_{\mathcal{Z}^{\lambda',\mu';1}_{t}}+\int^{t}_{0}\|w^{0}_{\hbar }(x-v\tau,v)\|_{\mathcal{Z}^{\lambda',\mu';1}_{\tau}}\|n_{\hbar}(\tau,\cdot)\|_{\mathcal{F}^{\lambda'\tau+\mu'}}d\tau.$$
 $$$$
 \begin{center}
\item\subsection{The property of  the Wigner function in the nonlinear case}
\end{center}

Sketch the proof of Theorem 0.1:

The  Newton iteration:

First of all, we write  a classical Newton iteration :
$$\textmd{Let}\quad\quad\quad\quad\quad\quad\quad\quad\quad\quad\quad\quad\quad\quad\quad\quad\quad\quad\quad\quad\quad\quad\quad\quad\quad\quad$$
$$ w_{\hbar}^{0}=w_{\hbar}^{0}(v) \quad\textrm{be\quad given},$$
and $$w_{\hbar}^{n}=w_{\hbar}^{0}+h_{\hbar}^{1}+\ldots+h_{\hbar}^{n},$$
where
 \begin{align}
  \left\{\begin{array}{l}
      \partial_{t}h_{\hbar}^{1}+v\cdot\nabla_{x}h_{\hbar}^{1}+\Theta_{\hbar}[\phi^{1}](h_{\hbar}^{1})w_{\hbar}^{0}=0,\\
     h_{\hbar}^{1}(0,x,v) =w_{\hbar I}-w_{\hbar}^{0}, \\
   \end{array}\right.
           \end{align}
 and now we consider the Vlasov equation in step $n+1,$  for any $n\geq1,$
 \begin{align}
  \left\{\begin{array}{l}
      \partial_{t}h_{\hbar}^{n+1}+v\cdot\nabla_{x}h_{\hbar}^{n+1}+\Theta_{\hbar}[\phi^{n}](w_{\hbar}^{n})h_{\hbar}^{n+1} \\
      =-\Theta_{\hbar}[\phi^{n+1}](h_{\hbar}^{n+1})w_{\hbar}^{n}-\Theta_{\hbar}[\phi^{n}](h_{\hbar}^{n})h_{\hbar}^{n},\\
      h_{\hbar}^{n+1}(0,x,v) =0.\\
   \end{array}\right.
           \end{align}


            In order to estimate the trajectory behavior of the particles, we have to  estimate the transition  from the quantum case to  the classical case. Indeed, since we restrict space variable on the localized domain $\mathbb{T}^{3} $ and  consider the problem under the stability condition and  the analytic frame,  then the trajectories of the particles have no difference  on the influence of the electric field under some analytic norm $F_{\lambda,\mu}$ between the quantum case and the classical case . In other words, under the quantum stability condition, we can transfer the discrete form in quantum case into the continuous form that is the same to the classical case. In fact, note that $$(\Theta_{\hbar}[\phi^{n}](w_{\hbar}^{n})h_{\hbar}^{n+1})(t,x,v)$$
           $$=\frac{im}{(2\pi)^{3}}\int_{\mathbb{R}^{3}}\int_{\mathbb{R}^{3}}\frac{\phi^{n}(w_{\hbar}^{n})(t,x+\frac{\hbar}{2m}\eta)
           -\phi^{n}(w_{\hbar}^{n})(t,x-\frac{\hbar}{2m}\eta)}
           {\hbar}h^{n+1}_{\hbar}(t,x,v')e^{i(v-v')\cdot\eta}dv'd\eta,$$
           by the inductive hypothesis, using the mean-value theorem,  we can obtain
           $$(\Theta_{\hbar}[\phi^{n}](w_{\hbar}^{n})h_{\hbar}^{n+1})(t,x,v)=C(x)\nabla_{x}\phi^{n}(w_{\hbar}^{n})(t,x)\int_{\mathbb{R}^{3}}i\eta
           h^{n+1}_{\hbar}(t,x,\eta)e^{iv\cdot\eta}d\eta$$
           $$=C(x)\nabla_{x}\phi^{n}(w_{\hbar}^{n})(t,x)\cdot\nabla_{v}h_{\hbar}^{n+1},$$
           where $C(x)$ are the uniform bounded on $x.$

   Therefore, we rewrite (4.13) into the form
    \begin{align}
  \left\{\begin{array}{l}
      \partial_{t}h_{\hbar}^{n+1}+v\cdot\nabla_{x}h_{\hbar}^{n+1}+C(x)\nabla_{x}\phi^{n}(w_{\hbar}^{n})\cdot\nabla_{v}h_{\hbar}^{n+1} \\
      =-\Theta_{\hbar}[\phi^{n+1}](h_{\hbar}^{n+1})w_{\hbar}^{n}-\Theta_{\hbar}[\phi^{n}](h_{\hbar}^{n})h_{\hbar}^{n},\\
      h_{\hbar}^{n+1}(0,x,v) =0.\\
   \end{array}\right.
           \end{align}

Then the corresponding continued trajectories  can be  described as follows:
for any $(x,v)\in\mathbb{T}^{3}\times\mathbb{R}^{3},$
let $(X^{n}_{\hbar;t,s},V^{n}_{\hbar;t,s})$ as the solution of the following ordinary differential equations
  $$
   \left\{\begin{array}{l}
   \frac{d}{dt}X^{n+1}_{\hbar;t,s}(x,v)=V^{n+1}_{\hbar;t,s}(x,v),\\
 X^{n+1}_{\hbar;s,s}(x,v)=x,
  \end{array}\right.
           $$

         \begin{align}
   \left\{\begin{array}{l}
   \frac{d}{dt}V^{n+1}_{\hbar;t,s}(x,v)=C(x)\nabla_{x}\phi^{n}(t,X^{n}_{\hbar;t,s}(x,v)),\\
 V^{n+1}_{\hbar;s,s}(x,v)=v.
  \end{array}\right.
           \end{align}
            At the same time, we consider the corresponding  linear dynamics system as follows,
            \begin{align}
   \left\{\begin{array}{l}
  \frac{d}{dt}X^{0}_{\hbar;t,s}(x,v)=V^{0}_{\hbar;t,s}(x,v),\quad  \frac{d}{dt}V^{0}_{\hbar;t,s}(x,v)=0,\\
 X^{0}_{\hbar;s,s}(x,v)=x,\quad V^{0}_{\hbar;s,s}(x,v)=v.
  \end{array}\right.
           \end{align}
             It is easy to check that
             $$\Omega^{n}_{\hbar;t,s}-Id\triangleq(\delta X^{n}_{\hbar;t,s},\delta V^{n}_{\hbar;t,s})\circ (X^{0}_{\hbar;s,t},V^{0}_{\hbar;s,t})
             = ( X^{n}_{\hbar;t,s}\circ(X^{0}_{\hbar;s,t},V^{0}_{\hbar;s,t})-Id, V^{n}_{\hbar;t,s}\circ(X^{0}_{\hbar;s,t},V^{0}_{\hbar;s,t})-Id).$$
            Therefore, in order to estimate $( X^{n}_{\hbar;t,s}\circ(X^{0}_{\hbar;s,t},V^{0}_{\hbar;s,t})-Id, V^{n}_{\hbar;t,s}\circ(X^{0}_{\hbar;s,t},V^{0}_{\hbar;s,t})-Id),$
             we only need to study $(\delta X^{n}_{\hbar;t,s},\delta V^{n}_{\hbar;t,s})\circ (X^{0}_{\hbar;s,t},V^{0}_{\hbar;s,t}).$

             From Eqs.(4.15) and (4.16),
             $$
   \left\{\begin{array}{l}
   \frac{d}{dt}\delta X^{n+1}_{\hbar;s,t}(x,v)=\delta V^{n+1}_{\hbar;s,t}(x,v),\\
 \delta X^{n+1}_{\hbar;s,s}(x,v)=0,
  \end{array}\right.
          $$
            \begin{align}
   \left\{\begin{array}{l}
   \frac{d}{dt}\delta V^{n+1}_{\hbar;s,t}(x,v)=\nabla_{x}\phi^{n}(w_{\hbar}^{n})(t,X^{n}_{\hbar;s,t}(x,v)),\\
\delta V^{n}_{\hbar;s,s}(x,v)=0.
  \end{array}\right.
           \end{align}
At present, we have reduced the Wigner-Poisson system into the classical Vlasov-Poisson system, the corresponding mathematical theory is the same, so we will give the corresponding results to be used without the proof,  because these results are totally the same with  the classical case, here we omit.

           Integrating (4.14) in time and $ h_{\hbar}^{n+1}(0,x,v) =0,$ we get
            \begin{align}
            h_{\hbar}^{n+1}(t,X^{n}_{\hbar;t,0}(x,v),V^{n}_{\hbar;t,0}(x,v))=\int^{t}_{0}\Sigma_{\hbar}^{n+1}(s,X^{n}_{\hbar;s,0}(x,v),V^{n}_{\hbar;s,0}(x,v))ds,
            \end{align}
           where$$\Sigma_{\hbar}^{n+1}(t,x,v)=-\Theta_{\hbar}[\phi^{n+1}](w_{\hbar}^{n})-\Theta_{\hbar}[\phi^{n}](h_{\hbar}^{n})
         .$$

           By the definition of $(X^{n}_{\hbar,t,s}(x,v),V^{n}_{\hbar,t,s}(x,v)),$ we have
           $$h_{\hbar}^{n+1}(t,x,v)=\int^{t}_{0}\Sigma_{\hbar}^{n+1}(s,X^{n}_{\hbar;s,t}(x,v),V^{n}_{\hbar;s,t}(x,v))ds$$
           $$=\int^{t}_{0}\Sigma_{\hbar}^{n+1}(s,\delta X^{n}_{\hbar;s,t}(x,v)+ X^{0}_{\hbar;s,t}(x,v),\delta V^{n}_{\hbar;s,t}(x,v)+ V^{0}_{\hbar;s,t}(x,v))ds.$$

           Since the unknown $h_{\hbar}^{n+1}$ appears on both sides of (4.18), we hope to get a self-consistent estimate. For this, we have little choice but to
           integrate in $v$ and get an integral equation on $n[h_{\hbar}^{n+1}]=\int_{\mathbb{R}^{3}}h_{\hbar}^{n+1}dv,$ namely

 \begin{align}
           &n[h_{\hbar}^{n+1}](t,x)=\int^{t}_{0}\int_{\mathbb{R}^{3}}(\Sigma_{\hbar}^{n+1}\circ\Omega^{n}_{\hbar;s,t}(x,v))(s, X^{0}_{\hbar;s,t}(x,v),V^{0}_{\hbar;s,t}(x,v))dvds\notag\\
           &=\int^{t}_{0}\int_{\mathbb{R}^{3}}-\bigg[\mathcal{E}^{n+1}_{\hbar;s,t}\cdot G_{\hbar;s,t}^{n}-\mathcal{E}^{n}_{\hbar;s,t}\cdot H^{n}_{\hbar;s,t}\bigg](s, x-v(t-s),v)\notag\\
           &=I_{\hbar}^{n+1,n}+II_{\hbar}^{n,n},
           \end{align}
  where
  $$
  \left\{\begin{array}{l}
   \mathcal{E}^{n+1}_{\hbar;s,t}\cdot  G_{\hbar;s,t}^{n}=[\Theta_{\hbar}[\phi^{n+1}](w_{\hbar}^{n})]\circ\Omega^{n}_{\hbar;s,t}(x,v),\\
   \quad \mathcal{E}^{n}_{\hbar;s,t}\cdot H^{n}_{\hbar;s,t}=[\Theta_{\hbar}[\phi^{n}](h_{\hbar}^{n})]\circ\Omega^{n}_{\hbar;s,t}(x,v).\\
   \end{array}\right.
   $$

   In the following we only sketch the main idea and main steps in the proof, the detailed proof can be found in [,,], here we omit. In this paper, we use the inductive method.

   The first step for $n=1,$
 it is known that (4.12) is a linear Vlasov equation. From section 4.1,  the conclusions of Theorem 0.1 hold.

 Now for any $i\leq n,i\in\mathbb{N}_{0},$ we assume that the following estimates hold,
 \begin{align}
 &\sup_{t\geq0}\|n[h_{\hbar}^{i}](t,\cdot)\|_{\mathcal{F}^{\lambda_{i}t+\mu_{i}}}\leq\delta_{i},\notag\\
 &\sup_{0\leq s\leq t}\|h^{i}_{\hbar;s}\circ\Omega^{i-1}_{\hbar;t,s}\|
_{\mathcal{Z}_{s-\frac{bt}{1+b}}^{\lambda_{i}(1+b),\mu_{i};1}}\leq\delta_{i},\notag\\
\end{align}
then we have the following inequalities, denote $(\mathbf{E}^{n}):$
$$\sup_{t\geq 0}\|[\nabla_{x}\phi^{i}](h_{\hbar}^{i})(t,\cdot)\|_{\mathcal{F}^{\lambda_{i}t+\mu_{i}}}<\delta_{i},$$
$$\sup_{0\leq s\leq t}\|\nabla_{x}(h^{i}_{\hbar;s}\circ\Omega^{i-1}_{\hbar;t,s})\|
_{\mathcal{Z}_{s-\frac{bt}{1+b}}^{\lambda_{i}(1+b),\mu_{i};1}}\leq\delta_{i},\sup_{0\leq s\leq t}\|(\nabla_{x}h^{i}_{\hbar;s})\circ\Omega^{i-1}_{\hbar;t,s}\|
_{\mathcal{Z}_{s-\frac{bt}{1+b}}^{\lambda_{i}(1+b),\mu_{i};1}}\leq\delta_{i},$$
$$\|(\nabla_{v}+s\nabla_{x})(h^{i}_{\hbar;s}\circ\Omega^{i-1}_{\hbar;t,s})\|
_{\mathcal{Z}_{s-\frac{bt}{1+b}}^{\lambda_{i}(1+b),\mu_{i};1}}\leq\delta_{i},\|((\nabla_{v}+s\nabla_{x})h^{i}_{\hbar;s})\circ\Omega^{i-1}_{\hbar;t,s})\|
_{\mathcal{Z}_{s-\frac{bt}{1+b}}^{\lambda_{i}(1+b),\mu_{i};1}}\leq\delta_{i},$$
$$\sup_{0\leq s\leq t}\frac{1}{(1+s)^{2}}\|(\nabla\nabla h^{i}_{\hbar;s})\circ\Omega^{i-1}_{\hbar;t,s}\|
_{\mathcal{Z}_{s-\frac{bt}{1+b}}^{\lambda_{i}(1+b),\mu_{i};1}}\leq\delta_{i},$$
$$\sup_{0\leq s\leq t}(1+s)^{2}\|(\nabla_{v} h_{\hbar}^{i})\circ\Omega^{i-1}_{\hbar;t,s}-\nabla_{v}(h_{\hbar}^{i}\circ\Omega^{i-1}_{\hbar;t,s})\|
_{\mathcal{Z}_{s-\frac{bt}{1+b}}^{\lambda_{i}(1+b),\mu_{i};1}}\leq\delta_{i}.$$


\begin{center}
   \item\subsection{Local time iteration}
\end{center}

Before working out the core of the proof of Theorem 0.1, we shall show a short time estimate, which will play a role as an initial data layer for the Newton scheme.
The main tool in this section is given by the following lemma.
\begin{lem} Let $f$ be an analytic function, $\lambda(t)=\lambda-Kt$ and $\mu(t)=\mu-Kt, K>0,$ let $T>0$ be so small that $\lambda(t)>0,
\mu(t)>0$ for $0\leq t\leq T.$ Then for any $s\in[0,T]$ and any $p\geq1,$
$$\frac{d^{+}}{dt}\bigg|_{t=s}\|f\|_{\mathcal{Z}_{s}^{\lambda(t),\mu(t);p}}\leq-\frac{K}{(1+s)}\|\nabla f\|_{\mathcal{Z}_{s}^{\lambda(s),\mu(s);p}},$$
where $\frac{d^{+}}{dt}$ stands for the upper right derivative.
\end{lem}

\begin{prop} There exists some small constant $T>0,$ such that when  all  conditions of Theorem 0.1 hold, then
for any fixed $\eta,k,$  for all $t\in(0,T],$ $0<\lambda<\lambda_{0},$  we have
         \begin{align}
         &|\hat{w}_{\hbar}(t,k,\eta)-\hat{w}_{\hbar I}(k,\eta)|\leq C(C_{0},\delta_{0}) e^{-2\pi\lambda|\eta+kt|}e^{-2\pi\mu |k|},\notag\\
         &|\hat{n}_{\hbar}(t,k)-\hat{n}_{\hbar I}|\leq C(C_{0},\delta_{0}) e^{-2\pi(\lambda t+\mu)|k|},\quad |\hat{E}(t,k)|\leq C(C_{0},\delta_{0})e^{-2\pi\lambda|k|t}e^{-2\pi\mu|k| },
         \end{align}
        where $\rho_{0}=\int_{\mathbb{T}^{3}}\int_{\mathbb{R}^{3}}f_{0}(x,v)dvdx.$
\end{prop}

\textbf{Proof.}  The first stage of the iteration,namely, $h^{1}$ will be  considered in the following section $\S4.1.$ So we only need to care about the higher orders.
Recall that $w_{\hbar}^{k}=w_{\hbar}^{0}+h_{\hbar}^{1}+\ldots+h_{\hbar}^{k}.$ And we define
$$\lambda_{k}(t)=\lambda_{k}-2Kt\quad \textmd{and}\quad \mu_{k}(t)=\mu_{k}-Kt,$$
where $\{\lambda_{k}\}^{\infty}_{k=1}$ and $\{\mu_{k}\}^{\infty}_{k=1}$ are decreasing sequences of positive numbers.

We assume inductively that at stage $ n$ of the iteration, we have constructed $\{\lambda_{k}\}^{n}_{k=1},\{\mu_{k}\}^{n}_{k=1},\{\delta_{k}\}^{n}_{k=1}$ such that
$$\sup_{0\leq t\leq T}\|h_{\hbar}^{k}(t,\cdot)\|_{\dot{\mathcal{Z}}^{\lambda_{k}(t),\mu_{k}(t);1}_{t}}\leq\delta_{k},$$
for all $ 1\leq k\leq n,$  for some fixed $T>0.$

In the following we need to show that the induction hypothesis are satisfied at stage $n+1.$ For this, we have to construct $\lambda_{n+1},\mu_{n+1},\delta_{n+1}.$

Note that $h^{n+1}_{\hbar}=0,$ at $t=0.$
 For $n\geq1,$ now let us solve
 $$ \partial_{t}h_{\hbar}^{n+1}+v\cdot\nabla_{x}h_{\hbar}^{n+1} =\widetilde{\Sigma}_{\hbar}^{n+1},$$
 where $$\widetilde{\Sigma}_{\hbar}^{n+1}=-\Theta_{\hbar}[\phi^{n}](w_{\hbar}^{n})h_{\hbar}^{n+1}-\Theta_{\hbar}[\phi^{n+1}](h_{\hbar}^{n+1})w_{\hbar}^{n}
-\Theta_{\hbar}[\phi^{n}](h_{\hbar}^{n})h_{\hbar}^{n}.$$

 Hence, $$\|h_{\hbar}^{n+1}\|_{\mathcal{Z}_{t}^{\lambda_{n+1}(t),\mu_{n+1}(t);1}}\leq\int^{t}_{0}\|\overline{\Sigma}_{\tau}^{n+1}\circ S_{-(t-s)}^{0}\|
 _{\mathcal{Z}_{t}^{\lambda_{n+1}(t),\mu_{n+1}(t);1}}ds\leq\int^{t}_{0}\|\overline{\Sigma}_{s}^{n+1}\|
 _{\mathcal{Z}_{s}^{\lambda_{n+1}(t),\mu_{n+1}(t);1}}ds,$$
  where $$\overline{\Sigma}^{n+1}=-\Theta_{\hbar}[\phi^{n}](w_{\hbar}^{n})h_{\hbar}^{n+1}-\Theta_{\hbar}[\phi^{n+1}](h_{\hbar}^{n+1})w_{\hbar}^{n}
-\Theta_{\hbar}[\phi^{n}](h_{\hbar}^{n})h_{\hbar}^{n}.$$
 then by Lemma 4.4,
 $$\frac{d^{+}}{dt}\|h_{\hbar}^{n+1}\|_{\mathcal{Z}_{t}^{\lambda_{n+1}(t),\mu_{n+1}(t);1}}\leq-K\|\nabla_{x}h_{\hbar}^{n+1}\|_{\mathcal{Z}_{t}
 ^{\lambda_{n+1},\mu_{n+1};1}}
 -K\|\nabla_{v}h_{\hbar}^{n+1}\|_{\mathcal{Z}_{t}^{\lambda_{n+1},\mu_{n+1};1}}$$
 $$+\|\nabla_{x}[\phi^{n}(w_{\hbar}^{n})\|_{\mathcal{F}^{\lambda_{n+1}t+\mu_{n+1}}}
 \|\nabla_{v}h_{\hbar}^{n+1}\|_{\mathcal{Z}_{t}^{\lambda_{n+1},\mu_{n+1};1}}
+\|\nabla_{x}\phi^{n+1}(h_{\hbar}^{n+1})\|_{\mathcal{F}^{\lambda_{n+1}t+\mu_{n+1}}}
 \|\nabla_{v}w_{\hbar}^{n}\|_{\mathcal{Z}_{t}^{\lambda_{n+1},\mu_{n+1};1}}$$
 $$+\|\nabla_{x}\phi^{n}(h_{\hbar}^{n})\|_{\mathcal{F}^{\lambda_{n+1}t+\mu_{n+1}}}
 \|\nabla_{v}h_{\hbar}^{n}\|_{\mathcal{Z}_{t}^{\lambda_{n+1},\mu_{n+1};1}}ds.$$

 First, we easily get $\|\Theta_{\hbar}[\phi^{n}]\|_{\mathcal{F}^{\lambda_{n+1}t+\mu_{n+1}}}\leq C\|\nabla h_{\hbar}^{n}\|_{\mathcal{Z}_{t}^{\lambda_{n+1},\mu_{n+1};1}}.$
Moreover,
 $$ \|\nabla_{v}w_{\hbar}^{n}\|_{\mathcal{Z}_{t}^{\lambda_{n+1},\mu_{n+1};1}}\leq\sum^{n}_{i=1}\|\nabla_{v}h_{\hbar}^{i}\|_{\mathcal{Z}_{t}^{\lambda_{n+1},\mu_{n+1};1}}
 \leq C\sum^{n}_{i=1}\frac{\|h_{\hbar}^{i}\|_{\mathcal{Z}_{t}^{\lambda_{i+1},\mu_{i+1};1}}}{\min\{\lambda_{i}-\lambda_{n+1},\mu_{i}-\mu_{n+1}\}}.$$

 We gather the above estimates,
 $$\frac{d^{+}}{dt}\|h_{\hbar}^{n+1}\|_{\mathcal{Z}_{t}^{\lambda_{n+1}(t),\mu_{n+1}(t);1}}\leq \bigg( C\sum^{n}_{i=1}\frac{\delta_{i}}
 {\min\{\lambda_{i}-\lambda_{n+1},\mu_{i}-\mu_{n+1}\}}-K\bigg)\|\nabla h_{\hbar}^{n+1}\|_{\mathcal{Z}_{t}^{\lambda_{n+1},\mu_{n+1};1}}$$
 $$
 +\frac{\delta^{2}_{n}}{\min\{\lambda_{n}-\lambda_{n+1},\mu_{n}-\mu_{n+1}\}}.$$

We may choose $$\delta_{n+1}=\frac{\delta^{2}_{n}}{\min\{\lambda_{n}-\bar{\lambda}_{n+1},\mu_{n}-\bar{\mu}_{n+1}\}},$$
if \begin{align}
C\sum^{n}_{i=1}\frac{\delta_{i}}
 {\min\{\lambda_{i}-\lambda_{n+1},\mu_{i}-\mu_{n+1}\}}\leq K\end{align}
  holds.

 We choose $\lambda_{i}-\lambda_{i+1}=\mu_{i}-\mu_{i+1}=\frac{\Lambda}{i^{2}},$ where $\Lambda>0$ is arbitrarily small. Then for $i\leq n,\lambda_{i}-\lambda_{n+1}\geq
 \frac{\Lambda}{i^{2}},$ and $\delta_{n+1}\leq\delta^{2}_{n}n^{2}/\Lambda.$ Next we need to check that $\sum^{\infty}_{n=1}\delta_{n}n^{2}<\infty.$
 In fact, we choose $ K $ large enough and T small enough such that $\lambda_{0}-KT\geq \lambda_{\ast},\mu_{0}-KT\geq \mu_{\ast},$ and (3.9) holds, where
  $\lambda_{0}>\lambda_{\ast},\mu_{0}>\mu_{\ast}$ are fixed.

  If $\delta_{1}=\delta,$ then $\delta_{n}=n^{2}\frac{\delta^{2^{n}}}{\Lambda^{n}}(2^{2})^{2^{n-2}}(4^{2})^{2^{n-2}}\ldots((n-1)^{2})^{2}n^{2}.$
  To prove the sequence convergence for $\delta$ small enough, by induction that $\delta_{n}\leq z^{a^{n}},$ where $z$ small enough and $a\in(1,2).$ We
  claim that the conclusion holds for $n+1.$ Indeed, $\delta_{n+1}\leq\frac{z^{2a^{n}}}{\Lambda}n^{2}\leq z^{a^{n+1}}\frac{z^{(2-a)a^{n}}n^{2}}{\Lambda}.$ If $ z$ is
  so small that $z^{(2-a)a^{n}}\leq\frac{\Lambda}{n^{2}}$ for all $n\in\mathbb{N},$ then $\delta_{n+1}\leq z^{a^{n+1}},$ this concludes the local-time argument.

   \begin{center}
\item\subsection{Global time iteration }
\end{center}

Based on the estimates of the local-time iteration, without loss of generality, sometimes  we only consider the case  $s\geq\frac{bt}{1+b},$ where $b$
is small enough.

First, we give deflection estimates that compare the free evolution with the true evolution for the particles trajectories.
  \begin{prop}
           Assume  for any $i\in\mathbb{N},0<i\leq n, $
           $$\sup_{t\geq 0}\|\nabla_{x}\phi^{i}(h_{\hbar}^{i})(t,\cdot)\|_{\mathcal{F}^{\lambda_{i}t+\mu_{i}}}<\delta_{i}.$$

           And there exist  constants $\lambda_{\star}>0,\mu_{\star}>0$ such that $\lambda_{0}>\lambda'_{0}>\lambda_{1}
           >\lambda'_{1}>\ldots>\lambda_{i}>\lambda'_{i}>\ldots>\lambda_{\star},
           \mu_{0}>\mu_{1}>\mu'_{1}>\ldots>\mu_{i}>\mu'_{i}>\ldots>\mu_{\star}.$




          Then we have 
          $$\|\delta X^{n+1}_{\hbar;t,s}\circ (X^{0}_{\hbar;s,t},V^{0}_{\hbar;s,t})\|_{\mathcal{Z}_{s-\frac{bt}{1+b}}^{\lambda'_{n},\mu'_{n}}}
          \leq C\sum^{n}_{i=1}\delta_{i}  e^{-\pi(\lambda_{i}-\lambda'_{i})s}\min\bigg\{\frac{(t-s)^{2}}{2},
           \frac{1}{2\pi(\lambda_{i}-\lambda'_{i})^{2}}\bigg\},$$

            $$\|\delta V^{n+1}_{\hbar;t,s}\circ (X^{0}_{\hbar;s,t},V^{0}_{\hbar;s,t})\|_{\mathcal{Z}_{s-\frac{bt}{1+b}}^{\lambda'_{n},\mu'_{n}}}
            \leq C\sum^{n}_{i=1}\delta_{i} e^{-\pi(\lambda_{i}-\lambda'_{i})s}\min\bigg\{\frac{(t-s)}{2},
           \frac{1}{2\pi(\lambda_{i}-\lambda'_{i})}\bigg\},$$
           for $0<s<t,b=b(t,s)$ sufficiently small.
           \end{prop}
\begin{prop} Under the assumptions of Proposition 4.4, then
           $$\bigg\|\nabla\Omega^{n+1}X_{\hbar;t,s}-(Id,0)\bigg\|_{\mathcal{Z}_{s+\frac{bt}{1-b}}^{\lambda'_{n+1}(1-b),\mu'_{n+1}}}<\mathcal{C}_{1}^{n},
           \bigg\|\nabla\Omega^{n+1}V_{\hbar;t,s}-(0,Id)\bigg\|_{\mathcal{Z}_{s+\frac{bt}{1-b}}^{\lambda'_{n+1}(1-b),\mu'_{n+1}}}<\mathcal{C}_{1}^{n}+
           \mathcal{C}_{2}^{n},$$
           where $\mathcal{C}_{1}^{n}= C\sum^{n}_{i=1}\frac{ e^{-\pi(\lambda_{i}-\lambda'_{i})s}\delta_{i}}
           {2\pi(\lambda_{i}-\lambda'_{i})^{2}}\min\bigg\{\frac{(t-s)^{2}}{2},
           1\bigg\},\mathcal{C}_{2}^{n}= C\sum^{n}_{i=1}\frac{ e^{-\pi(\lambda_{i}-\lambda'_{j})s}\delta_{i}}{2\pi(\lambda_{i}-\lambda'_{i})}
           \min\bigg\{t-s,1\bigg\}.$
           \end{prop}
           \begin{prop} Under the assumptions of Proposition 4.4, then
           $$\bigg\|\Omega^{i}X_{\hbar;t,s}-\Omega^{n}X_{\hbar;t,s}\bigg\|_{\mathcal{Z}_{s-\frac{bt}{1+b}}^{\lambda'_{n}(1-b),\mu'_{n}}}
           < \mathcal{C}_{1}^{i,n},\bigg\|\Omega^{i}V_{\hbar;t,s}-\Omega^{n}V_{\hbar;t,s}\bigg\|_{\mathcal{Z}_{s-\frac{bt}{1+b}}^{\lambda'_{n}(1-b),\mu'_{n}}}
           <\mathcal{C}_{1}^{i,n}+ \mathcal{C}_{2}^{i,n},$$
           where $\mathcal{C}_{1}^{i,n}=C\sum^{n}_{j=i+1} \frac{ e^{-\pi(\lambda_{j}-\lambda'_{j})s}\delta_{j}}
           {2\pi(\lambda_{j}-\lambda'_{j})^{2}}\min\bigg\{\frac{(t-s)^{2}}{2},
           1\bigg\},\mathcal{C}_{2}^{i,n}=C\sum^{n}_{j=i+1}
           \frac{ e^{-\pi(\lambda_{j}-\lambda'_{j})s}\delta_{j}}{2\pi(\lambda_{j}-\lambda'_{j})}
           \min\bigg\{t-s,1\bigg\}.$
           \end{prop}
           \begin{rem} Note that $\mathcal{C}_{1}^{i,n},\mathcal{C}_{2}^{i,n}$  decay fast as $s\rightarrow\infty,i\rightarrow\infty,$ and uniformly in
           $n\geq i,$ since the sequence $\{\delta_{n}\}^{\infty}_{n=1}$ has fast convergence. Hence, if $r\in\mathbb{N}$  given, we shall have
           \begin{align}
           \mathcal{C}_{1}^{i,n}\leq\omega_{i,n}^{r,1},\quad \textmd{and}\quad \mathcal{C}_{2}^{i,n}\leq\omega_{i,n}^{r,2},\quad \textmd{all}\quad r\geq1,
           \end{align}
           with $\omega_{i,n}^{r,1}=C^{r}_{\omega}\sum^{n}_{j=i+1} \frac{ \delta_{j}}
           {2\pi(\lambda_{j}-\lambda'_{j})^{2+r}}\frac{\min\{\frac{(t-s)^{2}}{2},
           1\}}{(1+s)^{r}}$ and $\omega_{i,n}^{r,2}=C^{r}_{\omega}\sum^{n}_{j=i+1} \frac{ \delta_{j}}
           {2\pi(\lambda_{j}-\lambda'_{j})^{1+r}}\frac{\min\{\frac{(t-s)^{2}}{2},
           1\}}{(1+s)^{r}},$ for some absolute constant $C^{r}_{\omega}$ depending only on $ r.$
           \end{rem}
           \begin{prop} Under the assumptions of Proposition 4.4, then
           $$\bigg\|(\Omega^{i}_{\hbar;t,s})^{-1}\circ\Omega_{\hbar;t,s}^{n}-Id\bigg\|_{\mathcal{Z}_{s-\frac{bt}{1+b}}^{\lambda'_{n}(1-b),\mu'_{n}}}
            <\mathcal{C}_{1}^{i,n}+ \mathcal{C}_{2}^{i,n}.$$
          \end{prop}

            To give a self-consistent estimate, we have to control each term of Eq.(4.19): $I_{\hbar}^{n+1,n}, II_{\hbar}^{n,n}.$ And the most difficult term is $I,$
            because there is some resonance phenomena occurring in this term that makes the  propagated wave away from equilibrium.

\begin{cor}  Under the assumptions of Theorem 0.1, $\alpha_{n}=\lambda_{n}-\lambda'_{n}, \varepsilon_{n}\in(0,\alpha_{n}),$
recalling that $\hat{\rho}(t,0)=0,$ and that our conditions imply an upper bound on $c_{n}$ and $c^{n}_{0},$ Under the following conditions hold for all $n\geq1$:
\begin{itemize}
\item[(I)]   $2C^{1}_{\omega}\bigg(\sum^{n}_{i=1}\frac{\delta_{i}}{(2\pi(\lambda_{i}-\lambda'_{n}))^{3}}\bigg)\leq\min\bigg\{\frac{\lambda^{\ast}_{n}B}{6},
\frac{\mu_{0}-\mu'_{n}}{2}\bigg\};$
         \item[(II)] $\mathcal{C}_{1}^{n}+\mathcal{C}_{2}^{n}\leq\varepsilon;$
         \item[(III)]    for all $i\in\{1,\ldots,n-1\}$ and all $t\geq\tau,$
         $2(1+\tau)(1+B)(3\mathcal{C}_{1}^{i,n}+\mathcal{C}_{2}^{i,n})(\tau,t)\leq\max\{\lambda'_{i}-\lambda'_{n},\mu'_{i}-\mu'_{n}\};$
          \item[(IV)] for all $i\in\{1,\ldots,n\}$ and all $\tau\in[0,t],$
          $4(1+\tau)(\mathcal{C}_{1}^{i,n}+\mathcal{C}_{2}^{i,n})\leq\min\{\lambda_{i}-\lambda'_{n},\mu_{i}-\mu'_{n}\};$
          \item[(V)] $4C^{1}_{\omega}\sum^{n}_{i=1}\frac{\delta_{i}}{(2\pi(\lambda_{i}-\lambda'_{i}))^{3}
       }\leq\frac{\lambda'_{\infty}B}{3};$
       \item[(VI)] $C\sum^{n}_{i=1}\frac{\delta_{i}}{(\lambda_{i}-\lambda'_{i})^{3}}\leq\frac{\lambda^{\ast}B}{3}-\iota;$
       \item[(VII)] $\bigg(C^{4}_{\omega}\bigg(C'_{0}+\sum^{n}_{i=1}\delta_{i}+1\bigg)\bigg(\sum^{n}_{j=1}\frac{\delta_{j}}
{2\pi(\lambda_{j}-\lambda'_{j})^{6}}\bigg)\leq\frac{1}{8};$
\item[(VIII)] $C_{W}\sum^{n}_{i=1}\frac{\delta_{i}}{\sqrt{2\pi(\lambda_{i}-\lambda'_{n})}}\leq\frac{1}{4},\quad
\sum^{n}_{i=1}\frac{\delta_{i}}{\pi(\lambda_{i}-\lambda'_{n})}\leq\max\{\chi,\frac{1}{8}\}.$
\end{itemize}

Then we have the uniform control,
$$\|\rho[h^{n+1}](t,\cdot)\|
_{\mathcal{F}^{\lambda'_{n}t+\mu'_{n}}}\leq\frac{C\delta^{2}_{n}(1+c^{n}_{0})^{2}}{\sqrt{\varepsilon_{n}}(\lambda_{n}-\lambda'_{n})^{2}}
\bigg(1+\frac{1}{\alpha_{n}\varepsilon_{n}^{\frac{3}{2}}}\bigg)e^{CT^{2}_{n}},$$
where $T_{n}=C\bigg(\frac{1}{\alpha^{5}_{n}\varepsilon_{n}}\bigg)^{\frac{1}{\gamma-1}}.$
\end{cor}

$Proof.$ From Propositions 4.4-4.6, we know that
$$\int^{t}_{0}K^{n}_{0}(t,s)ds\leq C_{W}\sum^{n}_{i=1}\frac{\delta_{i}}{\pi(\lambda_{i}-\lambda'_{n})},\quad
\int^{\infty}_{s}K^{n}_{0}(t,s)ds\leq C_{W}\sum^{n}_{i=1}\frac{\delta_{i}}{\pi(\lambda_{i}-\lambda'_{n})},$$
$$\bigg(\int^{t}_{0}K^{n}_{0}(t,s)^{2}ds\bigg)^{\frac{1}{2}}\leq C_{W}\sum^{n}_{i=1}\frac{\delta_{i}}{\sqrt{2\pi(\lambda_{i}-\lambda'_{n})}}.$$
Here $\alpha_{n}=\pi\min\{(\mu_{n}-\mu'_{n}),(\lambda_{n}-\lambda'_{n})\},$ and assume $\alpha_{n}$ is smaller than $\bar{\alpha}(\gamma)$ in Theorem 3.4, and
that
$$\bigg(C^{4}_{\omega}\bigg(C'_{0}+\sum^{n}_{i=1}\delta_{i}+1\bigg)\bigg(\sum^{n}_{j=1}\frac{\delta_{j}}
{2\pi(\lambda_{j}-\lambda'_{j})^{6}}\bigg)\leq\frac{1}{8}\quad\quad$$
$$C_{W}\sum^{n}_{i=1}\frac{\delta_{i}}{\sqrt{2\pi(\lambda_{i}-\lambda'_{n})}}\leq\frac{1}{4},\quad
\sum^{n}_{i=1}\frac{\delta_{i}}{\pi(\lambda_{i}-\lambda'_{n})}\leq\max\bigg\{\chi,\frac{1}{8}\bigg\}.\quad$$
Applying Theorem 3.4, we can deduce that for any $\varepsilon_{n}\in(0,\alpha_{n})$ and $t\geq0,$
$$\|\rho[h^{n+1}](t,\cdot)\|
_{\mathcal{F}^{\lambda'_{n}t+\mu'_{n}}}$$
$$\leq\frac{C\delta^{2}_{n}(1+c^{n}_{0})^{2}}{\sqrt{\varepsilon_{n}}(\lambda_{n}-\lambda'_{n})^{2}}
\bigg(1+\frac{1}{\alpha_{n}\varepsilon_{n}^{\frac{3}{2}}}\bigg)e^{CT^{2}_{n}},$$
where $T_{n}=C\bigg(\frac{1}{\alpha^{5}_{n}\varepsilon_{n}}\bigg)^{\frac{1}{\gamma-1}}.$

   \begin{center}
\item\subsection{  The proof of main theorem}
\end{center}

$$$$
$\mathbf{Step }$ $\mathbf{2}.$
 Without loss of generality, we assume that  $\hbar$ is small enough,  up to slightly lowering $\lambda_{1},$ we may choose all parameters in such a way that
$\lambda_{k},\lambda'_{k}\rightarrow\lambda_{\infty}>\underline{\lambda}\quad \textmd{and}\quad\mu_{k},
\mu'_{k}\rightarrow\mu_{\infty}>\underline{\mu},\quad \textmd{as}
\quad k\rightarrow\infty;$
then we pick up $B>0$ such that
$\mu_{\infty}-\lambda_{\infty}(1+B)B\geq\mu'_{\infty}>\underline{\mu},$
and we let $b(t)=\frac{B}{1+t}.$
From the iteration, we have,  for all $k\geq2,$
\begin{align}
\sup_{0\leq s\leq t}\|h^{k}_{\hbar;s}\circ\Omega^{k-1}_{\hbar;t,s}\|_{\mathcal{Z}^{\lambda_{\infty}(1+b),\mu_{\infty};1}_{t-\frac{bt}{1+b}}}\leq\delta_{k},
\end{align}
where $\sum^{\infty}_{k=2}\delta_{k}\leq C\delta.$
Choosing $t=s$ in (4.24) yields
$\sup_{0\leq s\leq t}\|h^{k}_{\hbar;s}\|_{\mathcal{Z}^{\lambda_{\infty}(1+B),\mu_{\infty};1}_{t-\frac{Bt}{1+B+t}}}\leq\delta_{k}.$
This implies that
$\sup_{ t\geq0}\|h^{k}_{\hbar;t}\|_{\mathcal{Z}^{\lambda_{\infty}(1+B),\mu_{\infty}-\lambda_{\infty}(1+B)B;1}_{t}}\leq\delta_{k}.$
In particular, we have a uniform estimate on $h^{k}_{\hbar;t}$ in $\mathcal{Z}^{\lambda_{\infty},\mu'_{\infty};1}_{t}.$ Summing up over $ k $ yields for
$w_{\hbar}=w_{\hbar}^{0}+\sum^{\infty}_{k=1}h_{\hbar}^{k},$ the estimate
\begin{align}
\sup_{t\geq0}\|w_{\hbar}(t,\cdot)-w_{\hbar}^{0}\|_{\mathcal{Z}^{\lambda_{\infty},\mu'_{\infty};1}_{t}}\leq C\delta.
\end{align}
From (viiii) of Proposition 1.3, we can deduce from (4.25) that
\begin{align}
\sup_{t\geq0}\|w_{\hbar}(t,\cdot)-w_{\hbar}^{0}\|_{\mathcal{Y}^{\underline{\lambda},\underline{\mu}}_{t}}\leq C\delta.
\end{align}
Moreover, $n_{\hbar}=\int_{\mathbb{R}^{3}}w_{\hbar}dv$ satisfies similarly
$\sup_{t\geq0}\|n_{\hbar}(t,\cdot)\|_{\mathcal{F}^{\lambda_{\infty}t+\mu_{\infty}}}\leq C\delta.$
It follows that $|\hat{n}_{\hbar}(t,k)|\leq C\delta $
$e^{-2\pi\lambda_{\infty}|k|t}e^{-2\pi\mu_{\infty}|k|}$ for any $k\neq 0.$ On the one hand, by Sobolev embedding,
we deduce that for any $r\in\mathbb{N},$
$$\|n_{\hbar}(t,\cdot)-\langle n_{\hbar}\rangle\|_{C^{r}(\mathbb{T}^{3})}\leq C_{r}\delta e^{-2\pi\lambda't};$$
on the other hand, multiplying $\hat{n}_{\hbar}$ by the Fourier transform of $W,$  we see that
$P$ satisfies
\begin{align}
|\widehat{\phi}(t,k)|\leq C\delta e^{-2\pi\lambda'|k|t}e^{-2\pi\mu'|k|};
\end{align}
for some $\lambda_{0}>\lambda'>\underline{\lambda},\mu_{0}>\mu'>\underline{\mu}.$

Now, from (4.26), we have, for any $(k,\eta)\in\mathbb{Z}^{3}\times\mathbb{R}^{3}$ and any $t\geq0,$
\begin{align}
|\hat{w}_{\hbar}(t,k,\eta+kt)-\hat{w}_{\hbar}^{0}(\eta)|\leq C\delta e^{-2\pi\mu'|k|}e^{-2\pi\lambda'|\eta|},
\end{align}
this finishes the proof of Theorem 0.1.

\begin{center}
\item\subsection{  The stability of the Schr$\ddot{\textmd{o}}$dinger equation}
\end{center}

In this section we go back the Schr$\ddot{\textmd{o}}$dinger equations.
\begin{align}
\left\{\begin{array}{l}i\hbar\partial_{t}\psi_{\alpha}=(-\frac{\hbar^{2}}{2m}\Delta_{x}-e\phi(t,x))\psi_{\alpha},\\
 \phi=4\pi W(x)\ast (n_{\hbar}-\bar{n}_{\hbar}),\\
\psi_{\alpha}(0,x)=\psi_{\alpha I}(x),
\end{array}\right.
\end{align}
then
\begin{align}
\psi_{\alpha}(t,x)=\frac{\hbar^{\frac{3}{2}}}{t^{\frac{3}{2}}}e^{-i\frac{|x|^{2}}{\hbar t}}\psi_{\alpha I}(x)+e\int^{t}_{0}\int_{\mathbb{\mathbb{T}}^{3}}\frac{\hbar^{\frac{3}{2}}}{(t-s)^{\frac{3}{2}}} e^{-i\frac{|x-y|^{2}}{\hbar (t-s)}}\phi(s,y)\psi_{\alpha}(s,y)dyds.
\end{align}

From Theorem 0.1, we have for any $0<\lambda''<\lambda',$
$$\|\psi_{\alpha}(t,\cdot)\|_{\mathcal{F}^{\lambda'' t+\mu}}$$
\begin{align}
&\leq\bigg\|\frac{\hbar^{\frac{3}{2}}}{t^{\frac{3}{2}}}e^{-i\frac{|x|^{2}}{\hbar t}}\psi_{\alpha I}\bigg\|_{\mathcal{F}^{\lambda'' t+\mu}}+e\int^{t}_{0}\bigg\|\frac{\hbar^{\frac{3}{2}}}{(t-s)^{\frac{3}{2}}} e^{-\frac{i|x|^{2}}{\hbar (t-s)}}\bigg\|_{\mathcal{F}^{\lambda s+\mu;\infty}}\|\phi(s,\cdot)\|_{\mathcal{F}^{\lambda'' s+\mu}}\|\psi_{\alpha}(s,\cdot)\|_{\mathcal{F}^{\lambda'' s+\mu}}ds\notag\\
&\leq \bigg\|\frac{\hbar^{\frac{3}{2}}}{t^{\frac{3}{2}}}e^{-i\frac{|x|^{2}}{\hbar t}}\psi_{\alpha I}\bigg\|_{\mathcal{F}^{\lambda'' t+\mu}}+C\int^{t}_{0}\bigg\|\frac{\hbar^{\frac{3}{2}}}{(t-s)^{\frac{3}{2}}} e^{-\frac{i|x|^{2}}{\hbar (t-s)}}\bigg\|_{\mathcal{F}^{\lambda s+\mu;\infty}}e^{-(\lambda'-\lambda'')s}\|\psi_{\alpha}(s,\cdot)\|_{\mathcal{F}^{\lambda'' s+\mu}}ds,\notag\\
\end{align}
then by Grownwell inequality,   \begin{align}
\|\psi_{\alpha}(t,\cdot)\|_{\mathcal{F}^{\lambda'' t+\mu}}\leq\bigg\|\frac{\hbar^{\frac{3}{2}}}{t^{\frac{3}{2}}}e^{-i\frac{|x|^{2}}{\hbar t}}\psi_{\alpha I}\bigg\|_{\mathcal{F}^{\lambda'' t+\mu}}\exp\bigg(\int^{t}_{0}\bigg\|\frac{\hbar^{\frac{3}{2}}}{(t-s)^{\frac{3}{2}}} e^{-\frac{i|x|^{2}}{\hbar (t-s)}}\bigg\|_{\mathcal{F}^{\lambda s+\mu;\infty}}e^{-(\lambda'-\lambda'')s}ds\bigg)
\end{align}

Acknowledgements: The author is grateful for the comfortable and superior academic environment of Yau Mathematical Science Center, Tsinghua University.  I also thank all including  my family and friends.Finally, I appreciate C.Villani and C.Mouhot's work, their work give a great help for me.

\end{document}